\documentclass[A4j,11pt]{article}
\usepackage{amsfonts,amsmath,amscd}
\usepackage{graphics}
\usepackage{amssymb}
\begin{document}

\title{{\bf A modfied mean curvature flow \\
in Euclidean space and soap bubbles \\
in symmetric spaces
}}
\author{{\bf Naoyuki Koike}}
\date{}
\maketitle

\begin{abstract}
In this paper, we show that small spherical soap bubbles in irreducible simply connected symmetric spaces of 
rank greater than one are constructed from the limits of a certain kind of modified mean curvature flows 
starting from small spheres in the Euclidean space of dimension equal to the rank of the symmetric space, 
where we note that the small spherical soap bubbles are invariant under the isotropy subgroup action of 
the isometry group of the symmetric space.  Furthermore, we investigate the shape and the mean curvature 
of the small spherical soap bubbles.  
\end{abstract}

\vspace{0.5truecm}

%\vspace{0.2truecm}

%\noindent
%$\overline{\qquad\qquad\qquad\qquad\qquad\qquad}$

%\vspace{0.1truecm}

%{\rm \footnotesize{2000 Mathematics Subject Classification. Primary 53C40, 
%53C35.}}

%{\rm \footnotesize{Key words and phrases. equifocal submanifold, polar 
%action.}}

\section{Introduction} 
Let $f$ be an immersion of an $n$-dimensional compact oriented manifold $M$ into 
an $(n+1)$-dimensional oriented Riemannian manifold $(\widetilde M,\widetilde g)$.  
If $f$ is of constant mean curvature, then $f(M)$ is called a {\it soap bubble}.  
Soap bubbles in the Euclidean space have studied by many geometers.  
In 1989, W.T. Hsiang and W.Y. Hsiang ([HH]) studied isoperimetric soap bubbles in the product space 
$H^{n_1}(c_1)\times\cdots\times H^{n_k}(c_k)$ of hyperbolic spaces or 
$H^{n_1}(c_1)\times\cdots\times H^{n_k}(c_k)\times{\Bbb R}^{n_{k+1}}$, where ``isoperimetric'' means that 
the soap bubble is a solution of the isopermetric problem.  
They proved that all isoperimetric soap bubbles in the product spaces are invariant under 
the isotropy subgroup action of the isometry group of the product space.  
Also, they ([HH]) proved that isoperimetric soap bubbles with the same constant mean curvature in 
$H^{n_1}(c_1)\times{\Bbb R}$ are congruent.  
Furthermore, they found the lower bound of the constant mean curvatures of isoperimetric soap bubbles in 
$H^{n_1}(c_1)\times H^{n_2}(c_2)$ or $H^{n_1}(c_1)\times{\Bbb R}^{n_2}$.  
In 1992, W.Y. Hsiang ([Hs]) found the lower bound 
%(which we denote by $b(G/K)$) 
of the constant mean curvatures of (not necessarily isoperimetric) soap bubbles in a rank ${\it l}(\geq 2)$ 
symmetric space $G/K$ of non-compact type and, in the case where $G/K$ is irreducible, he gave the explicit 
description of the lower bound by using the root system of $G/K$.  
Furthermore, he proved that, for each real number $b$ greater than the lower bound, there exists a $K$-invariant 
spherical soap bubble of constant mean curvature $b$ in $G/K$.  

In this paper, we introduce the notion of a weighted root system and define the modified 
mean curvature flow in a Euclidean space associated to the system.  
We show that the flows starting from small spheres exist in infinite time and converge to 
an embedded hypersurfaces in $C^{\infty}$-topology.  Furthermore, in the case where the system is 
one associated to an irreducible simply connected symmetric space of rank greater than one, 
we show that spherical soap bubbles in the symmetric spaces are constructed from the limit hypersurfaces of 
the flows starting from small spheres and investigate the shape and the mean curvature of the spherical soap 
bubbles.  

First we shall introduce the notion of a weighted root system.  
Let ${\cal S}=(V,(\triangle,$\newline
$\{m_{\alpha}\,\vert\,\alpha\in\triangle\},\varepsilon))$ be a system 
satisfying the following four conditions:

\vspace{0.2truecm}

(i) $V$ is a ${\it l}(\geq 2)$-dimensional real vector space equipped with an inner product 

$\langle\,\,,\,\,\rangle$,

(ii) $\triangle$ is a subset of the dual space $V^{\ast}$ of $V$ and it is a root system of 
rank ${\it l}$ in 

the sense of [He] (i.e., it is of type ${\mathfrak a}_{\it l}, {\mathfrak b}_{\it l}, {\mathfrak c}_{\it l}, 
{\mathfrak d}_{\it l}, \mathfrak{bd}_{\it l},
{\mathfrak e}_6$ (in case of ${\it l}=6$), ${\mathfrak e}_7$ (in 

case of ${\it l}=7$), ${\mathfrak e}_8$ (in case of ${\it l}=8$), ${\mathfrak k}_4$ (in case of ${\it l}=4$), 
${\mathfrak g}_2$ (in case of ${\it l}=2$)),

(iii) $m_{\alpha}$ ($\alpha\in\triangle)$ are positive integers and $\varepsilon$ is equal to $1$ or $-1$.  

\vspace{0.2truecm}

\noindent
We call the system ${\cal S}$ a {\it weighted root system} and ${\it l}$ the {\it rank} of 
the system.  Denote  by ${\rm rank}\,{\cal S}$ the rank of ${\cal S}$.  
Let $W$ be the group generated by the reflections with respect to $\alpha^{-1}(0)$'s ($\alpha\in\triangle$).  
We call $W$ the {\it reflection group associated with} ${\cal S}$.  
Let ${\cal S}_i=(V_i,(\triangle_i,\{m_{\alpha}\,\vert\,\alpha\in\triangle_i\},\varepsilon_i))$ ($i=1,2$) 
be weighted root systems.  
If $\varepsilon_1=\varepsilon_2$ and if there exists a linear isometry $\Phi$ of $V_1$ onto $V_2$ satisfying 

\vspace{0.2truecm}

(i) $\{\alpha\circ \Phi\,\vert\,\alpha\in\triangle_2\}=\triangle_1$, 

(ii) $m_{\alpha\circ \Phi}=m_{\alpha}$ ($\alpha\in\triangle_2$), 

\vspace{0.2truecm}

\noindent
then we say that ${\cal S}_1$ {\it is isomorphic to} ${\cal S}_2$ and call $\Phi$ an {\it isomorphism of} 
${\cal S}_1$ {\it onto} ${\cal S}_2$.  
%Since $\Phi(\Pi_1)=\Pi_2$, there exists a diffeomorphism $\overline{\Phi}$ of 
%$V_1/\Pi_1$ onto $V_2/\Pi_2$ with $\overline{\Phi}\circ\pi_1=\pi_2\circ\Phi$, where $\pi_i$ is the quotient map 
%of $V_i$ onto $V_i/\Pi_i$ ($i=1,2$).  
Let ${\cal S}=(V,(\triangle,\{m_{\alpha}\,\vert\,\alpha\in\triangle\},\varepsilon))$ be 
a weighted root systems.  
Let $S_V(r)$ is the sphere of radius $r$ centered at the origin in $V$.  
The reflection group $W$ preserves $S_V(r)$ invariantly and hence it acts on $S_V(r)$.  
Let $\triangle_+(\subset\triangle)$ be the positive root system under a lexicographic ordering of $V^{\ast}$ 
and $\delta(\in\triangle_+)$ the highest root.  See [He] about the definitions of the positive root system 
and the highest root.  
Define a positive number $r_{\cal S}$ by 
$$r_{\cal S}:=\left\{
\begin{array}{cc}
\displaystyle{\frac{\pi}{\vert\vert\delta\vert\vert}} & 
({\rm in}\,\,{\rm case}\,\,{\rm of}\,\,\varepsilon=1)\\
\displaystyle{\infty} & 
({\rm in}\,\,{\rm case}\,\,{\rm of}\,\,\varepsilon=-1),
\end{array}\right.$$
where $\vert\vert\delta\vert\vert$ is 
the norm of $\delta$ with respect to the inner product of $V^{\ast}$ induced from $\langle\,\,,\,\,\rangle$.  
Let $B_V(r_{\cal S})$ be the open ball of radius $r_{\cal S}$ centered at the origin in $V$ and 
set $B_{\cal S}:=B_V(r_{\cal S})\setminus\{0\}$.  
Fix $r\in(0,r_{\cal S})$.  
Denote by $C^{\infty}_W(S_V(r),V)$ (resp. $C^{\infty}_W(S_V(r),B_{\cal S})$) 
the space of all $W$-equivariant $C^{\infty}$-maps from $S_V(r)$ to $V$ (resp. $B_{\cal S}$) and 
by ${\rm Imm}^{\infty}_W(S_V(r),V)$ (resp. ${\rm Imm}^{\infty}_W(S_V(r),B_{\cal S})$) the space of all 
$W$-equivariant $C^{\infty}$-immersions of $S_V(r)$ into $V$ (resp. $B_{\cal S}$).  
For $\phi\in{\rm Imm}^{\infty}_W(S_V(r),B_{\cal S})$, define a $W$-invariant function $\rho_{\cal S,\phi}$ over 
$S_V(r)$ by 
$$\rho_{\cal S,\phi}(Z):=\sum_{\alpha\in\triangle_+}
\frac{m_{\alpha}\sqrt{\varepsilon}\alpha(\nu(Z))}
{\tan(\sqrt{\varepsilon}\alpha(\phi(Z))}\quad(Z\in S_V(r)),\leqno{(1.1)}$$
%define a function $\rho_{{\cal S},\phi}$ over $S_V(r)$ by 
%$\rho_{{\cal S},\phi}:=\rho_{\cal S}\circ(\phi\times\nu)\circ\delta_{S_V(r)}$, 
where $\nu(:S_V(r)\to S_V(1))$ is the Gauss map of $\phi$ (defined by assigning the outward unit normal vector of 
$\phi$ at $Z$ to each point $Z$ of $S_V(r)$).  
%, where $\delta_{S_V(r)}$ is the diagonal map for $S_V(r)$.  
Note that, if $\alpha(\phi(Z))=0$, then we have $\nu(Z)=\frac{\phi(Z)}{\vert\vert\phi(Z)\vert\vert}$ and hence 
$\frac{\sqrt{\varepsilon}\alpha(\nu(Z))}{\tan(\sqrt{\varepsilon}\alpha(\phi(Z)))}$ implies 
$\frac{1}{\vert\vert\phi(Z)\vert\vert}$.  
Define a map $D_{\cal S}$ from ${\rm Imm}^{\infty}_W(S_V(r),B_{\cal S})$ to $C^{\infty}_W(S_V(r),V)$ by 
$$D_{\cal S}(\phi)
:=\left(\frac{\int_{S_V(r)}\left(\vert\vert\Delta_{\widehat g}\phi\vert\vert+\rho_{\cal S,\phi}\right)
dv_{\widehat g}}{\int_{S_V(r)}dv_{\widehat g}}-(\vert\vert\Delta_{\widehat g}\phi\vert\vert+\rho_{{\cal S},\phi})
\right)\nu
\leqno{(1.2)}$$
%$$D^2_{\cal S}(\phi,\nu):=d\phi\left({\rm grad}_g(\vert\vert\Delta_g\phi\vert\vert
%+\rho_{\cal S}(\phi,\nu))\right)
%\leqno{(1.3)}$$
for $\phi\in{\rm Imm}^{\infty}_W(S_V(r),B_{\cal S})$, where $g$ is 
the $W$-invariant metric on $S_V(r)$ induced from $\langle\,\,,\,\,\rangle$ by $\phi$, $dv_{\widehat g}$ is the 
volume element of $\widehat g$, $\Delta_{\widehat g}$ is the Laplace operator with respect to $\widehat g$ and 
$\vert\vert\cdot\vert\vert$ is the norm of $(\cdot)$ with respect to $\langle\,\,,\,\,\rangle$.  
We consider the following evolution equation 
$$\frac{\partial\phi_t}{\partial t}=D_{\cal S}(\phi_t)\leqno{({\rm E}_{\cal S})}$$
for a $C^{\infty}$-family $\phi_t$ in ${\rm Imm}^{\infty}_W(S_V(r),B_{\cal S})$.  
Since $\vert\vert\Delta_{\widehat g_t}\phi_t\vert\vert$ is the mean curvature of $\phi_t$, this evolution 
equation is interpreted as a modified volume-preserving mean curvature flow equation in $V$, where 
$\widehat g_t$ is the $W$-invariant metric on $S_V(r)$ induced from $\langle\,\,,\,\,\rangle$ by $\phi_t$.  
Denote by $\iota_r$ the inclusion map of $S_V(r)$ into $V$.  
It is clear that $\iota_r$ is $W$-equivariant.  

First we prove the following result for the evolution equation $(E_{\cal S})$.  

\vspace{0.5truecm}

\noindent
{\bf Theorem A.} {\sl Let ${\cal S}=(V,(\triangle,\{m_{\alpha}\,\vert\,\alpha\in\triangle\},\varepsilon))$ be 
a weighted root system and $\iota_r$ the inclusion map of $S_V(r)$ into $V$.  
Then there exists a positive constant $R_0$ smaller than $\frac{r_{\cal S}}{8}$ 
such that, if $r<R_0$, 
then the solution $\phi_t$ of the evolution equation $({\rm E}_{\cal S})$ satisfying the initial condition 
$\phi_0=\iota_r$ uniquely exists in infinite time and $\phi_{t_i}$ 
%($t\in[0,\infty)$) remain to be strictly convex and 
converges to a $W$-equivariant $C^{\infty}$-embedding $\phi_{\infty}$ of $S_V(r)$ into 
$B_{\cal S}$ (in the $C^{\infty}$-topology) as $i\to\infty$ for some sequence $\{t_i\}_{i=1}^{\infty}$ in 
$[0,\infty)$ with $\lim_{i\to\infty}t_i=\infty$.  
Furthermore, $\phi_t$ ($0\leq t<\infty$) remain to be strictly convex and hence so is also $\phi_{\infty}$.}

\vspace{0.5truecm}

Let $G/K$ be an irreducible simply connected rank ${\it l}(\geq 2)$ symmetric space of compact type or 
non-compact type.  A weighted root system of rank ${\it l}$ is defined for $G/K$ in a natural manner 
(see Section 3).  We call this system the {\it weighted root system associated with} $G/K$ and denote it by 
${\cal S}_{G/K}$.  
Let ${\cal S}_{G/K}=(V,(\triangle,\{m_{\alpha}\,\vert\,\alpha\in\triangle\},\varepsilon))$, where we note that 
$$\varepsilon=\left\{
\begin{tabular}{ll}
$1$ & (when $G/K$ is of compact type)\\
$-1$ & (when $G/K$ is of non-compact type).\\
\end{tabular}
\right.$$
The vector space $V$ is identified with the tangent space $T_{eK}{\cal T}$ of a maximal flat totally geodesic 
submanifold ${\cal T}$ in $G/K$ through $eK$, where $e$ is the identity element of $G$.  
In the case where $G/K$ is of compact type, ${\cal T}$ is identified with 
the quotient space $V/\Pi$ of $V$ by a lattice $\Pi$ in $V$, and 
in the case where $G/K$ is of non-compact type, it is identified with $V$ oneself.  
For convenience, let $\Pi=\{{\bf 0}\}$ in the case where $G/K$ is of non-compact type.  
Let $W$ be the reflection group associated with ${\cal S}_{G/K}$.  
Then it is shown that $\Pi$ is $W$-invariant.  
Denote by $\pi$ the quotient map of $V$ onto $V/\Pi={\cal T}$ and 
$r(G/K)$ the injective radius of $G/K$.  
Note that $r(G/K)=\infty$ in the case where $G/K$ is of non-compact type.  
It is easy to show that $r(G/K)=r_{{\cal S}_{G/K}}$.  

The main theorem in this paper is as follows.  

\vspace{0.5truecm}

\noindent
{\bf Theorem B.} {\sl Let ${\cal S}=(V,(\triangle,\{m_{\alpha}\,\vert\,\alpha\in\triangle\},\varepsilon))$ 
and $\iota_r$ as in Theorem A.  
Assume that $r<R_0$ and ${\cal S}$ is isomorphic to one associated to an irreducible simply connected rank 
${\it l}(\geq 2)$ symmetric space $G/K$ of compact type or non-compact type, where $R_0$ is as in Theorem A.  
%Let $\phi_t$ ($t\in[0,\infty)$) be the solution of the evolution equation $({\rm E}_{\cal S})$ satisfying the 
%initial condition $\phi_0=\iota_r$ and $\phi_{\infty}$ the limit of the flow.  
Then the following statements $({\rm i})-({\rm iii})$ hold.  

{\rm (i)} The hypersurface $M:=K\cdot\pi(\phi_{\infty}(S_V(r)))$ in $G/K$ is a $K$-invariant strictly convex 
spherical soap bubble in $G/K$, where $\phi_{\infty}$ is as in Theorem A.  
%, where $\pi$ is the quotient map of $V$ onto $V/\Pi$ 
%and $V/\Pi$ is identified with a maximal flat totally geodesic submanifold in $G/K$ through $eK$.  

{\rm (ii)} Let $C(\subset V)$ be a Weyl domain (i.e., a fundamental domain of the reflection group $W$) 
and $\theta_0$ the element of $(0,\frac{\pi}{2})$ defined by 
$$\theta_0:=\mathop{\max}_{P}\mathop{\max}_{Z_1\in P_{\max}}\mathop{\max}_{Z_2\in P_{\min,Z_1}}
\,\angle\,Z_1{\bf 0}Z_2
\left(=\mathop{\max}_{P}\mathop{\max}_{Z_1\in P_{\min}}\mathop{\max}_{Z_2\in P_{\max,Z_1}}
\,\angle\,Z_1{\bf 0}Z_2\right),$$
where $\angle\,Z_1{\bf 0}Z_2$ denotes the angle between $\overrightarrow{{\bf 0}Z_1}$ and 
$\overrightarrow{{\bf 0}Z_2}$, 
$P$ moves over the Grassmannian of all two-planes in $V$, $P_{\max}$ (resp. $P_{\min}$) 
denotes the set of all the maximal (resp. minimal) points of the function 
$\rho_{{\cal S},\iota_r}$ over $S_V(r)\cap\overline C\cap P$ ($\overline C\,:\,$ the closure of $C$) 
and $P_{\max,Z_1}(\subset P_{\max})$ (resp. $P_{\min,Z_1}(\subset P_{\min})$) 
denotes the set (which is at most two-points set) of all the maximal (resp. minimal) points neighboring 
$Z_1$ in $S_V(r)\cap\overline C\cap P$ of the function.  
Then we have 
$$M\subset B\left(\frac{r}{\cos\theta_0}\right)\setminus B(r\cos\theta_0),$$
where $B(r\cos\theta_0)$ (resp. $B(\frac{r}{\cos\theta_0})$) is the closed geodesic ball of radius 
$r\cos\theta_0$ (resp. $\frac{r}{\cos\theta_0}$) in $G/K$ centered at $eK$.  

{\rm (iii)} Let $\eta_{\max}$ and $\eta_{\min}$ be the functions 
defined by 
$$\eta_{\max}(s):=\max_{Z\in S_V(1)}\left(\rho_{{\cal S},\iota_s}(Z)+\frac{{\it l}-1}{s}\right)$$
and 
$$\eta_{\min}(s):=\min_{Z\in S_V(1)}\left(\rho_{{\cal S},\iota_s}(Z)+\frac{{\it l}-1}{s}\right).$$
Then the constant mean curvature $H_M$ of $M$ satisfies 
$\eta_{\min}(r)\leq H_M\leq\eta_{\max}(r)$.}

\vspace{0.5truecm}

\noindent
{\it Remark 1.1.} For convenience, denote by $M(r)$ and $H(r)$ the soap bubble $M$ and the mean curvature $H_M$ 
as in the statement of Theorem A, respectively.  
Since the volume of the domain surrounded by $M(r)$ is strictly increasing (and continuous) with respect to $r$, 
it is shown that $H(r)$ is strictly decreasing (and continuous) with respect to $r$.  
Easily we can show $\lim\limits_{r\to 0}\eta_{\max}(r)=\lim\limits_{r\to 0}\eta_{\min}(r)=\infty$ and hence 
$\lim\limits_{r\to 0}H(r)=\infty$.  
On the other hand, in the case where $G/K$ is of non-compact type, we have 
$$\lim_{r\to\infty}\eta_{\max}(r)=\mathop{\max}_{Z\in S_V(1)}
\sum_{\alpha\in\triangle_+}m_{\alpha}\vert\alpha(Z)\vert(=:b_{\max}(G/K))$$
and
$$\lim_{r\to\infty}\eta_{\min}(r)=\mathop{\min}_{Z\in S_V(1)}
\sum_{\alpha\in\triangle_+}m_{\alpha}\vert\alpha(Z)\vert(=:b_{\min}(G/K)).$$
Hence we obtain 
$$b_{\min}(G/K)\leq\lim_{r\to\infty}H(r)\leq b_{\max}(G/K).$$

\vspace{0.5truecm}

By using (ii) of Theorem B, we can derive the following result.  

\vspace{0.5truecm}

\noindent
{\bf Corollary C.} {\sl Under the hypothesis of Theorem B, 
set 

\vspace{-0.3truecm}

$$\theta_{G/K}:=\mathop{\max}_{(Z_1,Z_2)\in S_V(r)\cap\overline C}\,\angle\,Z_1{\bf 0}Z_2.$$

\vspace{-0.3truecm}
Then we have 

\vspace{-0.3truecm}

$$M\subset B\left(\frac{r}{\cos\theta_{G/K}}\right)\setminus B(r\cos\theta_{G/K}).$$
}

\vspace{0.5truecm}

In particular, we obtain the following result in the case where $G/K$ is of rank two.  

\vspace{0.5truecm}

\noindent
{\bf Corollary D.} {\sl Under the hypothesis of Theorem B, 
assume that the rank of $G/K$ is equal to two.  
Then we have 
$$M\subset\left\{
\begin{array}{ll}
\displaystyle{B\left(\frac{2r}{\sqrt 3}\right)\setminus B\left(\frac{\sqrt 3r}{2}\right)} & 
(\triangle:(\mathfrak a_2){\rm-type}\,\,{\rm or}\,\,(\mathfrak g_2){\rm-type})\\
\displaystyle{B(\sqrt 2r)\setminus B\left(\frac{r}{\sqrt 2}\right)} & 
(\triangle:(\mathfrak b_2){\rm-type}).
\end{array}\right.
$$
%If $\triangle$ is of ($\mathfrak a_2$)-type or ($\mathfrak g_2$)-type, then we have 
%$M\subset B(r_1)\setminus B(\frac{\sqrt 3}{2}r_1)$ for some $r_1\in(r,\frac{2}{\sqrt 3}r)$.  
%Also, if $\triangle$ is of ($\mathfrak b_2$)-type, then we have 
%$M\subset B(r_1)\setminus B\left(\frac{r_1}{\sqrt 2}\right)$ for some $r_1\in(r,\sqrt 2r)$.
}

\vspace{0.5truecm}

According to Corollary D, when the root system of $G/K$ is of type $\mathfrak a_2$, 
$\phi_{\infty}(S_V(r))$ is as in Figure 1 for example.  

\vspace{0.5truecm}

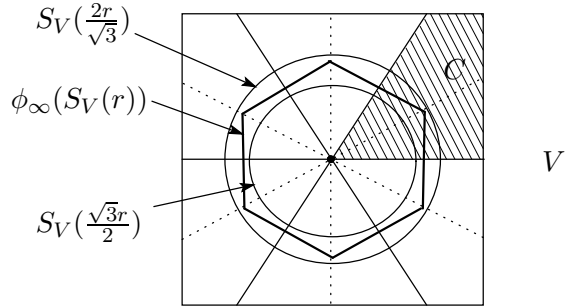
\centerline{
%WinTpicVersion3.08
\unitlength 0.1in
\begin{picture}( 44.0000, 18.7000)( -3.1000,-30.6000)
% BOX 2 0 3 0
% 2 2206 1200 3781 2725
% 
\special{pn 8}%
\special{pa 2206 1200}%
\special{pa 3782 1200}%
\special{pa 3782 2726}%
\special{pa 2206 2726}%
\special{pa 2206 1200}%
\special{fp}%
% DOT 0 0 3 0
% 2 2983 1962 2983 1962
% 
\special{pn 20}%
\special{sh 1}%
\special{ar 2984 1962 10 10 0  6.28318530717959E+0000}%
\special{sh 1}%
\special{ar 2984 1962 10 10 0  6.28318530717959E+0000}%
% LINE 2 2 3 0
% 2 3781 1962 3781 1962
% 
\special{pn 8}%
\special{pa 3782 1962}%
\special{pa 3782 1962}%
\special{dt 0.045}%
% LINE 2 0 3 0
% 2 2206 1962 3781 1962
% 
\special{pn 8}%
\special{pa 2206 1962}%
\special{pa 3782 1962}%
\special{fp}%
% LINE 2 0 3 0
% 2 2491 2725 3486 1200
% 
\special{pn 8}%
\special{pa 2492 2726}%
\special{pa 3486 1200}%
\special{fp}%
% LINE 2 0 3 0
% 2 3486 2725 2491 1200
% 
\special{pn 8}%
\special{pa 3486 2726}%
\special{pa 2492 1200}%
\special{fp}%
% LINE 3 0 3 0
% 54 3013 1962 3004 1945 3057 1962 3029 1908 3102 1962 3054 1870 3146 1962 3079 1832 3190 1962 3104 1794 3234 1962 3128 1757 3279 1962 3154 1719 3323 1962 3178 1681 3367 1962 3203 1644 3412 1962 3228 1606 3456 1962 3252 1568 3500 1962 3278 1530 3545 1962 3302 1493 3589 1962 3327 1455 3633 1962 3352 1418 3678 1962 3377 1380 3722 1962 3402 1342 3766 1962 3426 1304 3781 1905 3451 1266 3781 1819 3476 1229 3781 1733 3505 1200 3781 1648 3550 1200 3781 1562 3594 1200 3781 1476 3638 1200 3781 1390 3682 1200 3781 1304 3727 1200 3781 1219 3771 1200
% 
\special{pn 4}%
\special{pa 3014 1962}%
\special{pa 3004 1946}%
\special{fp}%
\special{pa 3058 1962}%
\special{pa 3030 1908}%
\special{fp}%
\special{pa 3102 1962}%
\special{pa 3054 1870}%
\special{fp}%
\special{pa 3146 1962}%
\special{pa 3080 1832}%
\special{fp}%
\special{pa 3190 1962}%
\special{pa 3104 1794}%
\special{fp}%
\special{pa 3234 1962}%
\special{pa 3128 1758}%
\special{fp}%
\special{pa 3280 1962}%
\special{pa 3154 1720}%
\special{fp}%
\special{pa 3324 1962}%
\special{pa 3178 1682}%
\special{fp}%
\special{pa 3368 1962}%
\special{pa 3204 1644}%
\special{fp}%
\special{pa 3412 1962}%
\special{pa 3228 1606}%
\special{fp}%
\special{pa 3456 1962}%
\special{pa 3252 1568}%
\special{fp}%
\special{pa 3500 1962}%
\special{pa 3278 1530}%
\special{fp}%
\special{pa 3546 1962}%
\special{pa 3302 1494}%
\special{fp}%
\special{pa 3590 1962}%
\special{pa 3328 1456}%
\special{fp}%
\special{pa 3634 1962}%
\special{pa 3352 1418}%
\special{fp}%
\special{pa 3678 1962}%
\special{pa 3378 1380}%
\special{fp}%
\special{pa 3722 1962}%
\special{pa 3402 1342}%
\special{fp}%
\special{pa 3766 1962}%
\special{pa 3426 1304}%
\special{fp}%
\special{pa 3782 1906}%
\special{pa 3452 1266}%
\special{fp}%
\special{pa 3782 1820}%
\special{pa 3476 1230}%
\special{fp}%
\special{pa 3782 1734}%
\special{pa 3506 1200}%
\special{fp}%
\special{pa 3782 1648}%
\special{pa 3550 1200}%
\special{fp}%
\special{pa 3782 1562}%
\special{pa 3594 1200}%
\special{fp}%
\special{pa 3782 1476}%
\special{pa 3638 1200}%
\special{fp}%
\special{pa 3782 1390}%
\special{pa 3682 1200}%
\special{fp}%
\special{pa 3782 1304}%
\special{pa 3728 1200}%
\special{fp}%
\special{pa 3782 1220}%
\special{pa 3772 1200}%
\special{fp}%
% STR 2 0 3 0
% 3 4090 1825 4090 1920 1 0
% $V$
\put(40.9000,-19.2000){\makebox(0,0)[lt]{$V$}}%
% VECTOR 2 0 3 0
% 2 2048 1657 2521 1829
% 
\special{pn 8}%
\special{pa 2048 1658}%
\special{pa 2522 1830}%
\special{fp}%
\special{sh 1}%
\special{pa 2522 1830}%
\special{pa 2466 1788}%
\special{pa 2472 1812}%
\special{pa 2452 1826}%
\special{pa 2522 1830}%
\special{fp}%
% STR 2 0 3 0
% 3 1940 1266 1940 1362 3 0
% $S_V(\frac{2r}{\sqrt 3})$
\put(19.4000,-13.6200){\makebox(0,0)[rb]{$S_V(\frac{2r}{\sqrt 3})$}}%
% STR 2 0 3 0
% 3 2028 1638 2028 1733 3 0
% $\phi_{\infty}(S_V(r))$
\put(20.2800,-17.3300){\makebox(0,0)[rb]{$\phi_{\infty}(S_V(r))$}}%
% STR 2 0 3 0
% 3 3574 1342 3574 1438 1 0
% $C$
\put(35.7400,-14.3800){\makebox(0,0)[lt]{$C$}}%
% LINE 2 2 3 0
% 2 2983 2715 2983 1190
% 
\special{pn 8}%
\special{pa 2984 2716}%
\special{pa 2984 1190}%
\special{dt 0.045}%
% LINE 2 2 3 0
% 2 2206 2381 3764 1543
% 
\special{pn 8}%
\special{pa 2206 2382}%
\special{pa 3764 1544}%
\special{dt 0.045}%
% LINE 2 2 3 0
% 2 3781 2372 2206 1562
% 
\special{pn 8}%
\special{pa 3782 2372}%
\special{pa 2206 1562}%
\special{dt 0.045}%
% LINE 1 0 3 0
% 2 3476 1714 2974 1447
% 
\special{pn 13}%
\special{pa 3476 1714}%
\special{pa 2974 1448}%
\special{fp}%
% LINE 1 0 3 0
% 2 2983 1447 2521 1724
% 
\special{pn 13}%
\special{pa 2984 1448}%
\special{pa 2522 1724}%
\special{fp}%
% LINE 1 0 3 0
% 2 2521 1724 2531 2219
% 
\special{pn 13}%
\special{pa 2522 1724}%
\special{pa 2532 2220}%
\special{fp}%
% LINE 1 0 3 0
% 4 2531 2219 2993 2477 2993 2477 3456 2219
% 
\special{pn 13}%
\special{pa 2532 2220}%
\special{pa 2994 2478}%
\special{fp}%
\special{pa 2994 2478}%
\special{pa 3456 2220}%
\special{fp}%
% LINE 1 0 3 0
% 2 3466 2219 3476 1714
% 
\special{pn 13}%
\special{pa 3466 2220}%
\special{pa 3476 1714}%
\special{fp}%
% ELLIPSE 2 0 3 0
% 4 2993 1972 3428 2368 3720 1972 4159 1972
% 
\special{pn 8}%
\special{ar 2994 1972 436 396  0.0000000 6.2831853}%
% ELLIPSE 2 0 3 0
% 4 2993 1962 3556 2507 4303 1962 4953 1962
% 
\special{pn 8}%
\special{ar 2994 1962 564 546  0.0000000 6.2831853}%
% STR 2 0 3 0
% 3 1999 2086 1999 2181 4 0
% $S_V(\frac{\sqrt 3r}{2})$
\put(19.9900,-21.8100){\makebox(0,0)[rt]{$S_V(\frac{\sqrt 3r}{2})$}}%
% VECTOR 2 0 3 0
% 2 2058 2210 2580 2095
% 
\special{pn 8}%
\special{pa 2058 2210}%
\special{pa 2580 2096}%
\special{fp}%
\special{sh 1}%
\special{pa 2580 2096}%
\special{pa 2512 2090}%
\special{pa 2528 2106}%
\special{pa 2520 2130}%
\special{pa 2580 2096}%
\special{fp}%
% STR 2 0 3 0
% 3 1860 2760 1860 2860 1 0
% {\scriptsize The six corners of $\phi_{\infty}(S_V(r))$ are smooth.}
\put(18.6000,-28.6000){\makebox(0,0)[lt]{{\scriptsize The six corners of $\phi_{\infty}(S_V(r))$ are smooth.}}}%
% VECTOR 2 0 3 0
% 2 1959 1304 2590 1581
% 
\special{pn 8}%
\special{pa 1960 1304}%
\special{pa 2590 1582}%
\special{fp}%
\special{sh 1}%
\special{pa 2590 1582}%
\special{pa 2538 1536}%
\special{pa 2542 1560}%
\special{pa 2522 1574}%
\special{pa 2590 1582}%
\special{fp}%
% STR 2 0 3 0
% 3 1860 2960 1860 3060 1 0
% {\scriptsize The six edges of $\phi_{\infty}(S_V(r))$ are curved.}
\put(18.6000,-30.6000){\makebox(0,0)[lt]{{\scriptsize The six edges of $\phi_{\infty}(S_V(r))$ are curved.}}}%
\end{picture}%
\hspace{7truecm}}

\vspace{0.5truecm}

\centerline{{\bf Figure 1.}}

\section{The volume-preserving mean curvature flow}
In this section, we shall recall the definition of the volume-preserving mean curvature flow and 
the result of N.D. Alikakos and A. Freire ([AF]) for this flow.  
Let $M$ be an $n$-dimensional compact oriented manifold, $(\widetilde M,\widetilde g)$ be 
an $(n+1)$-dimensional oriented Riemannian manifold and 
$f$ an immersion of $M$ into $\widetilde M$.  
Also, let $\{f_t\}_{t\in[0,T)}$ be a $C^{\infty}$-family of 
immersions of $M$ into $\widetilde M$, 
% and $N_t$ be the unit normal vector field of $f_t$ compatible 
%with the orietations of $M$ and $\widetilde M$
where $T$ is a positive constant or $T=\infty$.  Define a map 
$F:M\times[0,T)\to\widetilde M$ by $F(x,t)=f_t(x)$ 
($(x,t)\in M\times[0,T)$).  
Denote by $\pi_M$ the natural projection of $M\times[0,T)$ onto $M$.  
For a vector bundle $E$ over $M$, denote by $\pi_M^{\ast}E$ the induced bundle of $E$ 
by $\pi_M$.  
Denote by $g_t$ and $N_t$ the induced metric and the outward unit normal vector of $f_t$, 
respectively.  Also, denote by $H_t$ the mean curvature of $f_t$ for $-N_t$.  
Define the function $H$ over $M\times[0,T)$ by $H_{(x,t)}:=(H_t)_x$ ($(x,t)\in M\times [0,T)$), 
the section $g$ of $\pi_M^{\ast}(T^{(0,2)}M)$ by 
$g_{(x,t)}:=(g_t)_x$ ($(x,t)\in M\times [0,T)$) and the section $N$ of 
$\pi_M^{\ast}(T\widetilde M)$ by $N_{(x,t)}:=(N_t)_x$ ($(x,t)\in M\times [0,T)$), where 
$T^{(0,2)}M$ is the tensor bundle of $(0,2)$-type of $M$ and $T\widetilde M$ is the tangent 
bundle of $\widetilde M$.  
The average mean curvature $\overline H(:[0,T)\to{\Bbb R})$ is defined by 
$$\overline H_t:=\frac{\int_MH_tdv_{g_t}}{\int_Mdv_{g_t}},$$
where $dv_{g_t}$ is the volume element of $g_t$.  
G. Huisken ([Hu3]) called the flow $\{f_t\}_{t\in[0,T)}$ a {\it volume-preserving mean curvature flow} if 
it satisfies 
$$F_{\ast}\left(\frac{\partial}{\partial t}\right)=(\overline H-H)N.$$
He studied this flow in [Hu3].  
Along this flow, the volume of $(M,g_t)$ decreases but 
the volume of the domain $D_t$ sorrounded by $M_t:=f_t(M)$ is preserved invariantly.  
In particular, if $f_t$'s are embeddings, then we call $\{M_t\}_{t\in[0,T)}$ 
rather than $\{f_t\}_{t\in[0,T)}$ a volume-preserving mean curvature flow.  

Assume that $\widetilde M$ is compact.  Let $S$ be a geodesic sphere in $\widetilde M$ such that, 
for any point $p$ of the domain surrounded by $S$, $S$ lies in a geodesically convex domain of $p$ 
(in $\widetilde M$), and $\iota_S$ the inclusion map of $S$ into $\widetilde M$.  
Then, according to Main Theorem of [AF], it is shown that, if $S$ is a small geodesic sphere of radius 
smaller than some positive constant among the above geodesic spheres, then 
for any strictly convex $C^{\infty}$-embedding $f$ of $S$ into $\widetilde M$ which is sufficiently close to 
$\iota_S$, there exists the volume-preserving mean curvature flow $f_t$ starting from $f$ in infinite time, 
each $f_t$ is strictly convex and $f_{t_i}$ converges to a strictly convex $C^{\infty}$-embedding $f_{\infty}$ 
of constant mean curvature (in $C^{\infty}$-topology) as $i\to\infty$ for some sequence 
$\{t_i\}_{i=1}^{\infty}$ with $\lim_{i\to\infty}t_i=\infty$.  
Furthermore, if all critical points of the scalar curvature functions of $f_t$'s are non-degenerate, then 
$f_t$ converges to the embedding $f_{\infty}$ (in $C^{\infty}$-topology) as $t\to\infty$.  
Note that the positive constants $\delta_i$ and $\varepsilon_i$ ($i=0,1,2,3$) in the statement of Main Theorem of 
[AF] are taken as 
$$\delta_3<\delta_2<\delta_1/2<\delta_0/4\,\,\,\,{\rm and}\,\,\,\,\varepsilon_3<\varepsilon_2<\varepsilon_1/2
<\varepsilon_0/4$$
(see P258,291 and 296 of [AF]) and that $\delta_0$ and $\varepsilon_0$ are taken as in P257 of [AF].  
By the compactness of $\widetilde M$, the existenceness of the positive constant $\delta_0$ in Page 257 of [AF] 
is assured.  Hence the existencenesses of $\delta_i\,(i=1,2,3)$ also are assured.  
In the case where $\widetilde M$ is homogeneous, it is clear that their existencenesses are assured even if 
it is not compact.  Hence the statement of Main Theorem of [AF] is valid in the case where $\widetilde M$ is 
a (not necessarily compact) homogeneous space.  
Also, we note that the statement of Main Theorem of [AF] does not hold without the assumption that $f$ is 
sufficiently close to $\iota_S$.  
For example, in the case where $f$ is a strictly convex embedding of a small geodesic sphere $S$ into 
$\widetilde M=S^{n+1}(1)$ (which is not close to $\iota_S$) as in Figure 2, 
the volume-preserving mean curvature flow $f_t$ starting from $f$ ($0\leq t<\infty$) does not remain to be 
strictly convex.  This fact is stated in Remarks of Page 38 of [Hu3].  

\vspace{0.2truecm}

\centerline{
%WinTpicVersion3.08
\unitlength 0.1in
\begin{picture}( 14.5000, 13.1000)( 18.8000,-21.0000)
% CIRCLE 2 0 3 0
% 4 2592 1402 3102 1722 3432 1402 3802 1402
% 
\special{pn 8}%
\special{ar 2592 1402 602 602  0.0000000 6.2831853}%
% ELLIPSE 2 2 3 0
% 4 2592 1402 3192 1512 3412 1402 1822 1402
% 
\special{pn 8}%
\special{ar 2592 1402 600 110  3.1415927 3.1753955}%
\special{ar 2592 1402 600 110  3.2768039 3.3106067}%
\special{ar 2592 1402 600 110  3.4120152 3.4458180}%
\special{ar 2592 1402 600 110  3.5472265 3.5810293}%
\special{ar 2592 1402 600 110  3.6824377 3.7162405}%
\special{ar 2592 1402 600 110  3.8176490 3.8514518}%
\special{ar 2592 1402 600 110  3.9528603 3.9866631}%
\special{ar 2592 1402 600 110  4.0880715 4.1218743}%
\special{ar 2592 1402 600 110  4.2232828 4.2570856}%
\special{ar 2592 1402 600 110  4.3584941 4.3922969}%
\special{ar 2592 1402 600 110  4.4937053 4.5275081}%
\special{ar 2592 1402 600 110  4.6289166 4.6627194}%
\special{ar 2592 1402 600 110  4.7641279 4.7979307}%
\special{ar 2592 1402 600 110  4.8993391 4.9331419}%
\special{ar 2592 1402 600 110  5.0345504 5.0683532}%
\special{ar 2592 1402 600 110  5.1697617 5.2035645}%
\special{ar 2592 1402 600 110  5.3049729 5.3387758}%
\special{ar 2592 1402 600 110  5.4401842 5.4739870}%
\special{ar 2592 1402 600 110  5.5753955 5.6091983}%
\special{ar 2592 1402 600 110  5.7106067 5.7444096}%
\special{ar 2592 1402 600 110  5.8458180 5.8796208}%
\special{ar 2592 1402 600 110  5.9810293 6.0148321}%
\special{ar 2592 1402 600 110  6.1162405 6.1500434}%
\special{ar 2592 1402 600 110  6.2514518 6.2831853}%
% ELLIPSE 2 0 3 0
% 4 2592 1402 3192 1512 1832 1402 3582 1402
% 
\special{pn 8}%
\special{ar 2592 1402 600 110  6.2831853 6.2831853}%
\special{ar 2592 1402 600 110  0.0000000 3.1415927}%
% CIRCLE 1 0 3 0
% 4 2602 1412 2682 1682 3062 1412 2342 1412
% 
\special{pn 13}%
\special{ar 2602 1412 282 282  3.1415927 6.2831853}%
% ELLIPSE 1 0 3 0
% 4 2602 1412 2882 1492 2082 1412 3312 1412
% 
\special{pn 13}%
\special{ar 2602 1412 280 80  6.2831853 6.2831853}%
\special{ar 2602 1412 280 80  0.0000000 3.1415927}%
% ELLIPSE 1 0 3 0
% 4 2601 1338 2830 1512 3456 1557 1954 1524
% 
\special{pn 13}%
\special{ar 2602 1338 230 174  2.7796076 6.2831853}%
\special{ar 2602 1338 230 174  0.0000000 0.3249051}%
% ELLIPSE 1 0 3 0
% 4 2597 1402 2894 1550 2299 1665 2958 1729
% 
\special{pn 13}%
\special{ar 2598 1402 298 148  1.0677191 2.0846233}%
% STR 2 0 3 0
% 3 3330 860 3330 960 2 0
% $f(S)$ (strictly convex)
\put(33.3000,-9.6000){\makebox(0,0)[lb]{$f(S)$ (strictly convex)}}%
% STR 2 0 3 0
% 3 3320 1140 3320 1240 2 0
% $f_t(S)$ (not strictly convex)
\put(33.2000,-12.4000){\makebox(0,0)[lb]{$f_t(S)$ (not strictly convex)}}%
% STR 2 0 3 0
% 3 3220 1730 3220 1830 1 0
% $\widetilde M=S^{n+1}(1)$ 
\put(32.2000,-18.3000){\makebox(0,0)[lt]{$\widetilde M=S^{n+1}(1)$ }}%
% VECTOR 2 0 3 0
% 2 3290 1180 2830 1350
% 
\special{pn 8}%
\special{pa 3290 1180}%
\special{pa 2830 1350}%
\special{fp}%
\special{sh 1}%
\special{pa 2830 1350}%
\special{pa 2900 1346}%
\special{pa 2880 1332}%
\special{pa 2886 1308}%
\special{pa 2830 1350}%
\special{fp}%
% STR 2 0 3 0
% 3 1880 2000 1880 2100 1 0
% $f$ is not close to $\iota_S$.
\put(18.8000,-21.0000){\makebox(0,0)[lt]{$f$ is not close to $\iota_S$.}}%
% VECTOR 2 0 3 0
% 2 3300 910 2810 1210
% 
\special{pn 8}%
\special{pa 3300 910}%
\special{pa 2810 1210}%
\special{fp}%
\special{sh 1}%
\special{pa 2810 1210}%
\special{pa 2878 1192}%
\special{pa 2856 1182}%
\special{pa 2856 1158}%
\special{pa 2810 1210}%
\special{fp}%
% ELLIPSE 1 0 3 0
% 4 2610 1310 2850 1600 2230 1440 2350 1670
% 
\special{pn 13}%
\special{ar 2610 1310 240 290  2.2881991 2.8646838}%
% ELLIPSE 1 0 3 0
% 4 2590 1320 2350 1610 2850 1680 2970 1450
% 
\special{pn 13}%
\special{ar 2590 1320 240 290  0.2769088 0.8533935}%
\end{picture}%
\hspace{3truecm}}

\vspace{0.5truecm}

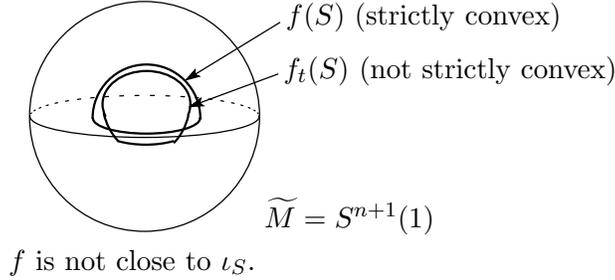
\centerline{{\bf Figure 2.}}

%\newpage\vspace{0.5truecm}

\section{The mean curvatures of hypersurfaces invariant under the isotropy action}
Let $G/K$ be an irreducible simply connected symmetric space of compact type or non-compact type.  
In this section, we shall first define the weighted root system associated with $G/K$.  
Set $n:={\rm dim}(G/K)-1$ and ${\it l}:={\rm rank}(G/K)$.  
Let $\mathfrak g$ (resp. $\mathfrak k$) be the Lie algebra of $G$ (resp. $K$) and 
$\mathfrak g=\mathfrak k+\mathfrak p$ the canonical decomposition associated with the symmetric pair 
$(G,K)$.  The space $\mathfrak p$ is identified with the tangent space $T_{eK}(G/K)$ of 
$G/K$ at $eK$, where $e$ is the identity element of $G$.  
Take a maximal abelian subspace $\mathfrak a$ of $\mathfrak p$.  
Let $\triangle(\subset\mathfrak a^{\ast})$ be the (restricted) root system of the symmetric pair 
$(G,K)$ with respect to $\mathfrak a$ and $\mathfrak p_{\alpha}$ the root space for 
$\alpha\in\triangle$.  See [He] about the definitions of these notions.  
Then we have the following root space decomposition:
$$\mathfrak p=\mathfrak a\oplus\left(\mathop{\oplus}_{\alpha\in\triangle_+}\mathfrak p_{\alpha}\right),$$
where $\triangle_+$ is the positive root system under some lexicographic ordering of 
$\mathfrak a^{\ast}$.  Set $m_{\alpha}:={\rm dim}\,\mathfrak p_{\alpha}$ ($\alpha\in\triangle$).  
Let $W$ be the Weyl group of $G/K$ with respect to $\mathfrak a$ (i.e., 
the group generated by the reflections with respect to the hyperplanes $\alpha^{-1}(0)$'s ($\alpha\in\triangle$) 
in $\mathfrak a$), which acts on 
$\mathfrak a$, $C(\subset\mathfrak a)$ be a Weyl domain (i.e., a fundamental domain of the action 
$W\curvearrowright\mathfrak a$) and $S_{\mathfrak a}(r)$ be 
the sphere of radius $r$ centered at the origin ${\bf 0}$ in $\mathfrak a$.  
The isotropy group $K$ acts on $G/K$ naturally.  This action is called the isotropy action of 
$G/K$.  Denote by $\exp$ the exponential map of $G$ and ${\rm Exp}$ the exponential map of 
$G/K$ at $eK$.  
Set ${\cal T}:={\rm Exp}(\mathfrak a)$, which is a maximal flat totally geodesic submanifold of $G/K$.  
Note that 
%${\rm dim}\,{\cal T}={\it l}$ and that ${\rm Exp}\vert_{\mathfrak a}:\mathfrak a\to{\cal T}$ is 
%the universal Riemannian covering and 
${\cal T}$ is identified with the quotient $\mathfrak a/\Pi$ by a $W$-invariant lattice $\Pi$ in $\mathfrak a$ 
in the case where $G/K$ is of compact type, and 
it is identified with $\mathfrak a$ oneself in the case where $G/K$ is of non-compact type.  
Set 
$$\varepsilon:=\left\{
\begin{tabular}{ll}
$1$ & (when $G/K$ is of compact type)\\
$-1$ & (when $G/K$ is of non-compact type).\\
\end{tabular}
\right.\leqno{(3.1)}$$
Then the system ${\cal S}:=(\mathfrak a,(\triangle,\{m_{\alpha}\,\vert\,\alpha\in\triangle\},\varepsilon))$ 
is a weighted root system.  We call this system the {\it weighted root system associated with} $G/K$ 
and denote it by ${\cal S}_{G/K}$.  

Next we shall describe explicitly the mean curvatures of $K$-invariant hypersurfaces in $G/K$ 
in terms of the roots.  
As a special case, those of geodesic spheres in $G/K$ are described explicitly.  
Let $M^{\mathfrak a}$ be a $W$-invariant star-shaped hypersurface (at the origin ${\bf 0}$) in $\mathfrak a$.  
Assume that, in the case where $G/K$ is of compact type, 
$\max_{Z\in M^{\mathfrak a}}d({\bf 0},Z)$ is smaller than the injective radius $r(G/K)$ of $G/K$, where 
$d$ is the Euclidean distance of $\mathfrak a$.  
Set $M_{\cal T}:={\rm Exp}(M^{\mathfrak a})$ and 
$M:=K\cdot M_{\cal T}(=K\cdot({\rm Exp}(M^{\mathfrak a}\cap\overline C))$.  Note that $M_{\cal T}$ and $M$ are hypersurfaces in ${\cal T}$ and $G/K$, respectively.  
Denote by $N$ the outward unit normal vector field of $M(\subset G/K)$, and 
$A$ and $H$ the shape operator and the mean curvature of $M(\subset G/K)$ 
for the inward unit normal vector field $-N$, respectively.  Also, denote by 
$A^{\cal T}$ and $H_{\cal T}$ those of $M_{\cal T}(\subset {\cal T})$ for 
the inward unit normal vector $-N\vert_{M_{\cal T}}$, respectively, and 
$\widehat A$ and $\widehat H$ those of $M^{\mathfrak a}(\subset\mathfrak a)$ for 
the inward unit normal vector, respectively.  
Take any $Z\in M^{\mathfrak a}\cap C$.  
Denote by $L_Z$ the $K$-orbit $K\cdot{\rm Exp}\,Z$, which is a principal orbit of the $K$-action 
because $Z\in C$.  We have 
$$T_{{\rm Exp}\,Z}L_Z=(\exp\,Z)_{\ast}
\left(\mathop{\oplus}_{\alpha\in\triangle_+}\mathfrak p_{\alpha}\right).$$
Denote by $A^Z$ the shape tensor of $L_Z(\subset G/K)$.  
Take $X_{\alpha}\in\mathfrak p_{\alpha}$ ($\alpha\in\triangle_+$).  
Then we have 
$$A^Z_{N_{{\rm Exp}\,Z}}((\exp\,Z)_{\ast}X_{\alpha})
=-\frac{\sqrt{\varepsilon}\alpha((\exp\,Z)_{\ast}^{-1}(N_{{\rm Exp}\,Z}))}
{\tan(\sqrt{\varepsilon}\alpha(Z))}(\exp\,Z)_{\ast}X_{\alpha}\leqno{(3.2)}$$
by suitably rescaling of the metric of $G/K$ (see [GT],[K1],[K2]), 
where $\frac{1}{\tan(\sqrt{\varepsilon}\alpha(Z))}=0$ in case of $\sqrt{\varepsilon}\alpha(Z)=\frac{\pi}{2}$.  
In the sequel, we give $G/K$ this rescaled metric.  
The vector field $N\vert_{L_Z}$ is a $K$-equivariant normal vector field along $L_Z$.  Hence, since 
the $K$-action is a hyperpolar action, it is a parallel normal vector field of $L_Z$ 
(see Theorem 5.5.12 of [PT]).  Hence, we have 
$$A^Z_{N_{{\rm Exp}(Z)}}=-A_{{\rm Exp}\,Z}\vert_{T_{{\rm Exp}\,Z}L_Z},\leqno{(3.3)}$$
where we note that $A_{{\rm Exp}\,Z}$ denotes the value of $A$ at ${\rm Exp}\,Z$.  In the sequel, 
we express $A_{{\rm Exp}\,Z}$ as $A$ for simplicity.  
From $(3.2)$ and $(3.3)$, we obtain 
$$A((\exp\,Z)_{\ast}X_{\alpha})
=\frac{\sqrt{\varepsilon}\alpha((\exp\,Z)_{\ast}^{-1}(N_{{\rm Exp}\,Z}))}
{\tan(\sqrt{\varepsilon}\alpha(Z))}(\exp\,Z)_{\ast}X_{\alpha}.\leqno{(3.4)}$$
Take $X_0\in\mathfrak a\ominus{\rm Span}\{N_{{\rm Exp}\,Z}\}$.  
Since ${\cal T}$ is totally geodesic in $G/K$, 
we have 
$$A((\exp\,Z)_{\ast}X_0)=A^{\cal T}((\exp\,Z)_{\ast}X_0)=(\exp\,Z)_{\ast}(\widehat AX_0).
\leqno{(3.5)}$$
From $(3.4)$ and $(3.5)$, we obtain 
$$H_{{\rm Exp}\,Z}=\sum_{\alpha\in\triangle_+}
\frac{m_{\alpha}\sqrt{\varepsilon}\alpha((\exp\,Z)_{\ast}^{-1}(N_{{\rm Exp}\,Z}))}
{\tan(\sqrt{\varepsilon}\alpha(Z))}+\widehat H_Z\qquad(Z\in M^{\mathfrak a}\cap C).\leqno{(3.6)}$$
Take any $Z'\in M^{\mathfrak a}\cap\partial C$, where $\partial C$ is the boundary of $C$.  
Set $\triangle^{Z'}_+:=\{\alpha\in\triangle_+\,\vert\,\alpha(Z')=0\}$.  
Then it follows from $(3.6)$ and the continuity of $H$ that 
$$H_{{\rm Exp}\,Z'}=\sum_{\alpha\in\triangle_+\setminus\triangle_+^{Z'}}
\frac{m_{\alpha}\sqrt{\varepsilon}\alpha((\exp\,Z')_{\ast}^{-1}(N_{{\rm Exp}\,Z'}))}
{\tan(\sqrt{\varepsilon}\alpha(Z'))}
+\frac{1}{\vert\vert Z'\vert\vert}\sum_{\alpha\in\triangle_+^{Z'}}m_{\alpha}+\widehat H_{Z'},\leqno{(3.7)}$$
where we note that $(\exp\,Z')_{\ast}^{-1}(N_{{\rm Exp}\,Z'})=\frac{1}{\vert\vert Z'\vert\vert}Z'$.  
The hypersuface $M$ is of constant mean curvature $\kappa$ if and only if 
$H_{{\rm Exp}\,Z}=\kappa$ ($Z\in M^{\mathfrak a}$) because $M$ is $K$-invariant.  
We consider the case where $G/K$ is of rank two.  Let $c:[0,b)\to\mathfrak a$ be the curve parametrized by 
the arclength whose image is equal to $M^{\mathfrak a}$, where $b$ is the length of $M^{\mathfrak a}$.  
Then it follows from $(3.6)$ that $M$ is of constant mean curvature $\kappa$ 
if and only if the following relation holds:
$$-\sum_{\alpha\in\triangle_+}
\frac{m_{\alpha}\sqrt{\varepsilon}\alpha(c''(s)/\vert\vert c''(s)\vert\vert)}
{\tan(\sqrt{\varepsilon}\alpha(c(s)))}+c''(s)=\kappa\quad\,\,(s\in[0,b)).$$
This relation is equivalent to the relation $(26)$ of [Hs, Page 164].  
W.Y. Hsiang ([Hs]) derived some facts by using the relation $(26)$.  

In particular, we consider the case where $M$ is a geodesic sphere.  
Let $S(r)$ be the geodesic sphere of radius $r(>0)$ centered at $eK$ and $S_{\mathfrak a}(r)$ 
the sphere of radius $r$ centered at the origin ${\bf 0}$ in $\mathfrak a$.  
Assume that $r$ is smaller than the first conjugate radius of $G/K$ in the case where 
$G/K$ is of compact type.  
%Set $S_{\cal T}(r):={\rm Exp}(S_{\mathfrak a}(r))$.  
Then we have $S(r)=K\cdot({\rm Exp}(S_{\mathfrak a}(r)\cap\overline C))$.  
Denote by $N^r$ the outward unit normal vector of $S(r)$ and $H^r$ the mean curvature of $S(r)$ for $-N^r$.  
Take any $Z\in S_{\mathfrak a}(r)\cap\overline C$.  
Since $(N^r)_{{\rm Exp}\,Z}=\frac{1}{r}(\exp\,Z)_{\ast}(Z)$, it follows from $(3.6)$ and $(3.7)$ that 
$$(H^r)_{{\rm Exp}\,Z}=\left\{
\begin{array}{ll}
\displaystyle{\sum_{\alpha\in\triangle_+}
\frac{m_{\alpha}\sqrt{\varepsilon}\alpha(Z)}{r\tan(\sqrt{\varepsilon}\alpha(Z))}
+\frac{{\it l}-1}{r}} & (Z\in S_{\mathfrak a}(r)\cap C)\\
\displaystyle{\begin{array}{l}
\displaystyle{\sum_{\alpha\in\triangle_+\setminus\triangle_+^{Z}}
\frac{m_{\alpha}\sqrt{\varepsilon}\alpha(Z)}
{r\tan(\sqrt{\varepsilon}\alpha(Z))}}\\
\displaystyle{+\frac{1}{r}\left(\sum_{\alpha\in\triangle_+^{Z}}m_{\alpha}+{\it l}-1\right)}
\end{array}} & 
(Z\in S_{\mathfrak a}(r)\cap\partial C).
\end{array}\right.\leqno{(3.8)}$$

\section{Proof of Theorem A}
In this section, we shall prove Theorem A stated in Introduction.  
We use the notations in Sections 1-3.  
Let ${\cal S}=(V,(\triangle,\{m_{\alpha}\,\vert\,\alpha\in\triangle\},\varepsilon)),\,\phi_t$ and $\iota_r$ 
be as in the statement of Theorem A.  
Assume that there exists a solution $\phi_t$ ($t\in[0,T)$) of the $(E_{\cal S})$ satisfying $\phi_0=\iota_r$.  
Denote by $\widehat g_t$ and $\nu_t$ the induced metric and the outward unit normal vector of $\phi_t$, 
respectively.  
Also, denote by $\widehat h_t$ and $\widehat A_t$ the second fundamental form and the shape operator of $\phi_t$ 
for $-\nu_t$, respectively.  Define the sections $\widehat g$ (resp. $\widehat h$) of 
$\pi_{S_V(r)}^{\ast}(T^{(0,2)}S_V(r))$ by 
$\widehat g_{(x,t)}:=(\widehat g_t)_x$ (resp. $\widehat h_{(x,t)}:=(\widehat h_t)_x$) 
($(x,t)\in S_V(r)\times [0,T)$), 
where $\pi_{S_V(r)}$ denotes the natural projection of $S_V(r)\times[0,T)$ onto $S_V(r)$.  
Also, denote by $\widehat{\nabla}^t$ the Riemannian connection of $\widehat g_t$ and 
$\Delta_{\widehat g_t}$ the Laplace operator with respect to $\widehat g_t$.  
Define a function $\widehat r_t$ ($t\in[0,T)$) over $S_V(r)$ by 
$$\widehat r_t(Z):=\vert\vert\phi_t(Z)\vert\vert\qquad(Z\in S_V(r))$$
and a diffeomorphism $\widehat c_t$ of $S_V(r)$ by 
$$\widehat c_t(Z):=\frac{r\phi_t(Z)}{\vert\vert\phi_t(Z)\vert\vert}\qquad(Z\in S_V(r)).$$
Then we can derive the following fact for the evolution of $\widehat r_t$.  

\vspace{0.5truecm}

\noindent
{\bf Lemma 4.1.} {\sl The functions $\{\widehat r_t\}_{t\in[0,T)}$ satisfies the following evolution equation:
$$\begin{array}{l}
\displaystyle{\frac{\partial\widehat r}{\partial t}
=\left(
\frac{\int_{S_V(r)}\left(\vert\vert\Delta_{\widehat g_t}\phi_t\vert\vert+\rho_{{\cal S},\phi_t}\right)
dv_{\widehat g_t}}{\int_{S_V(r)}dv_{\widehat g_t}}-(\vert\vert\Delta_{\widehat g_t}\phi_t\vert\vert
+\rho_{{\cal S},\phi_t})\right)}\\
\hspace{1.5truecm}\displaystyle{\times\frac{\widehat r_t\cdot\vert\vert c_{t\ast}
({\rm grad}_t\widehat r_t)\vert\vert}{\sqrt{\vert\vert{\rm grad}_t\widehat r_t\vert\vert^4+\widehat r_t^2
\vert\vert c_{t\ast}({\rm grad}_t\widehat r_t)\vert\vert^2}}.}
\end{array}\leqno{(4.1)}$$
}

\vspace{0.5truecm}

\noindent
{\it Proof.} By a simple calculation, we have 
$$\begin{array}{l}
\displaystyle{\frac{\partial\widehat r}{\partial t}=\langle D_{\cal S}(\phi_t),\frac{1}{\widehat r_t}\phi_t
\rangle}\\
\hspace{0.65truecm}\displaystyle{=\left(
\frac{\int_{S_V(r)}\left(\vert\vert\Delta_{\widehat g_t}\phi_t\vert\vert+\rho_{{\cal S},\phi_t}\right)
dv_{\widehat g_t}}{\int_{S_V(r)}dv_{\widehat g_t}}-(\vert\vert\Delta_{\widehat g_t}\phi_t\vert\vert
+\rho_{{\cal S},\phi_t})\right)}\\
\hspace{1.15truecm}\displaystyle{\times\frac{\langle\nu_t,\phi_t\rangle}{\widehat r_t}}
\end{array}$$
Also, by a simple calculation, we have 
$$\begin{array}{l}
\displaystyle{\nu_t=-\frac{\vert\vert{\rm grad}_t\widehat r_t\vert\vert^2}
{\vert\vert c_{t\ast}({\rm  grad}_t\widehat r_t)\vert\vert
\sqrt{\vert\vert{\rm grad}_t\widehat r_t\vert\vert^4+\widehat r_t^2
\vert\vert c_{t\ast}({\rm grad}_t\widehat r_t)\vert\vert^2}}\cdot
c_{t\ast}({\rm grad}_t\widehat r_t)}\\
\hspace{0.9truecm}\displaystyle{+\frac{\vert\vert c_{t\ast}({\rm grad}_t\widehat r_t)\vert\vert}
{\sqrt{\vert\vert{\rm grad}_t\widehat r_t\vert\vert^4+\widehat r_t^2
\vert\vert c_{t\ast}({\rm grad}_t\widehat r_t)\vert\vert^2}}\cdot\phi_t.}
\end{array}$$
From these relations, we obtain the desired evolution equation.  \hspace{2.65truecm}q.e.d.

\vspace{0.5truecm}

The following evolution equation holds for $\widehat h_t$.  

\vspace{0.5truecm}

\noindent
{\bf Lemma 4.2.} {\sl The families $\{\widehat h_t\}_{t\in[0,\infty)}$ satisfies 
$$\begin{array}{l}
\displaystyle{\frac{\partial\widehat h}{\partial t}-\Delta_{\widehat g_t}\widehat h_t
=\widehat{\nabla}^td(\rho_{{\cal S},\phi_t})
+{\rm Tr}(\widehat A_t^2)\widehat h_t}\\
\hspace{2.5truecm}\displaystyle{
+\left(\frac{\int_{S_V(r)}\left(\vert\vert\Delta_{\widehat g_t}\phi_t\vert\vert+\rho_{{\cal S},\phi_t}\right)
dv_{\widehat g_t}}{\int_{S_V(r)}dv_{\widehat g_t}}-2\vert\vert\Delta_{\widehat g_t}\phi_t\vert\vert
-\rho_{{\cal S},\phi_t}\right)\widehat h_t(\widehat A_t\bullet,\bullet).}
\end{array}$$
}

\vspace{0.5truecm}

\noindent
{\it Proof.} For simplicity, set $\widehat H^{\cal S}_t:=\vert\vert\Delta_{\widehat g_t}\phi_t\vert\vert
+\rho_{{\cal S},\phi_t}$ and 
$$\overline H^{\cal S}_t:=\frac{\int_{S_V(r)}\left(\vert\vert\Delta_{\widehat g_t}\phi_t\vert\vert
+\rho_{{\cal S},\phi_t}\right)dv_{\widehat g_t}}{\int_{S_V(r)}dv_{\widehat g_t}}.$$
By a simple calculation, we have 
$$\frac{\partial\widehat h}{\partial t}=\widehat{\nabla}^td\widehat H^{\cal S}_t
+(\overline H^{\cal S}_t-\widehat H^{\cal S}_t)\widehat h_t(\widehat A_t\bullet,\bullet).$$
Also, by using the Simon's identity, we have 
$$\Delta_{\widehat g_t}\widehat h_t=\widehat{\nabla}^td\widehat H_t+\widehat H_t
\widehat h_t(\widehat A_t\bullet,\bullet)-{\rm Tr}(\widehat A_t^2)\widehat h_t.$$
From these relations, we obtain the desired evolution equation.  \hspace{2.5truecm}q.e.d.

\vspace{0.5truecm}

Also, we prepare the following lemma, which will be used in the proof of the statement (iii) of Theorem B.  

\vspace{0.5truecm}

\noindent
{\bf Lemma 4.3.} {\sl Let $\overline H^{\cal S}_t$ and $\widehat H^{\cal S}_t$ be as in the proof of Lemma 4.2.  
The family $\{\widehat H^{\cal S}_t\}_{t\in[0,\infty)}$ satisfies the following evolution equation:
$$\begin{array}{l}
\displaystyle{\frac{\partial\widehat H^{\cal S}}{\partial t}-\Delta_{\widehat g_t}\widehat H^{\cal S}_t
=\sum_{\alpha\in\triangle_+}\frac{m_{\alpha}\sqrt{\varepsilon}\alpha(\phi_{t\ast}({\rm grad}_{\widehat g_t}
\widehat H_t^{\cal S}))}{\tan(\sqrt{\varepsilon}(\alpha\circ\phi_t))}}\\
\hspace{1.25truecm}\displaystyle{+(\overline H^{\cal S}_t-\widehat H^{\cal S}_t)
\sum_{\alpha\in\triangle_+}\frac{m_{\alpha}\sqrt{\varepsilon}^2(\alpha\circ\nu_t)^2
(3\cos^2(\sqrt{\varepsilon}(\alpha\circ\phi_t))-1)}
{\sin^2(\sqrt{\varepsilon}(\alpha\circ\phi_t))}.}
\end{array}\leqno{(4.2)}$$
}

\vspace{0.5truecm}

\noindent
{\it Proof.} The family $\{\widehat g_t\}_{t\in[0,\infty)}$ satisfies 
$$\frac{\partial\widehat g}{\partial t}=2(\overline H^{\cal S}_t-\widehat H^{\cal S}_t)\widehat h_t.$$
From the evolution equation in Lemma 4.2 and this evolution equation, we have 
$$\frac{\partial\widehat H^{\cal S}}{\partial t}-\Delta_{\widehat g_t}\widehat H^{\cal S}_t
=\frac{\partial\rho_{{\cal S},\phi_t}}{\partial t}+3(\overline H^{\cal S}_t-\widehat H^{\cal S}_t)
{\rm Tr}(\widehat A_t^2).$$
On the other hand, we have 
$$
%\begin{array}{l}\displaystyle{
\frac{\partial\rho_{{\cal S},\phi_t}}{\partial t}=
\sum_{\alpha\in\triangle_+}m_{\alpha}\left(-\frac{(\overline H^{\cal S}_t-\widehat H^{\cal S}_t)
\sqrt{\varepsilon}^2(\alpha\circ\nu_t)^2}{\sin^2(\sqrt{\varepsilon}(\alpha\circ\phi_t))}
+\frac{\sqrt{\varepsilon}\alpha(\phi_{t\ast}({\rm grad}_{\widehat g_t}\widehat H^{\cal S}_t))}
{\tan(\sqrt{\varepsilon}(\alpha\circ\phi_t))}\right).
%}\end{array}
$$
From these relations, we can derive the desired evolution equation.  
%$$\begin{array}{l}
%\displaystyle{\frac{\partial\widehat H^{\cal S}}{\partial t}-\Delta_{\widehat g_t}\widehat H^{\cal S}_t
%=(\overline H^{\cal S}_t-\widehat H^{\cal S}_t)
%\sum_{\alpha\in\triangle_+}m_{\alpha}\frac{\sqrt{\varepsilon}^2(\alpha\circ\nu_t)^2
%(3\cos^2(\sqrt{\varepsilon}(\alpha\circ\phi_t))-1)}
%{\sin^2(\sqrt{\varepsilon}(\alpha\circ\phi_t))}\\
%\hspace{2truecm}\displaystyle{
%+\sum_{\alpha\in\triangle_+}\frac{m_{\alpha}
%\sqrt{\varepsilon}\alpha(\phi_{t\ast}({rm grad}_{\widehat g_t}\widehat H^{\cal S}_t))}
%{\tan(\sqrt{\varepsilon}(\alpha\circ\phi_t))}.}
%\end{array}$$
\hspace{1.8truecm}q.e.d.

\vspace{0.5truecm}

By using Lemmas 4.1 and 4.2, we shall prove Theorem A.  

\vspace{0.5truecm}

\noindent
{\it Proof of Theorem A.} 
Since $L_t:=\phi_t(S_V(r))$ ($t\in[0,T)$) are $W$-invariant, their barycenter are equal to the origin ${\bf 0}$ 
of $V$.  Hence the barycenter $\xi(t)$ of $L_t$ is equal to the origin ${\bf 0}$ of $V$ and the diffeomorphisms 
$e(t)$ in the barysentric system $(3.1)$ in Page 288 of [AF] are regarded as 
the identity transformation of $S_V(1)$ under the identification of $T_{\bf 0}V$ and $V$.  
Hence the left-hand sides of the first and the second relations in $(3.1)$ are equal to zero.  
On the other hand, it is clear that the right-hand sides in the first and the second relations 
are equal to zero in our setting.  Thus the first and the second relations in $(3.1)$ (of [AF]) are trivial.  
Also, it is easy to show that $(4.1)$ corresponds to the third relation in $(3.1)$ (of [AF]), where we regard 
$\widehat r_t$ as a function over $S_V(1)$ under the natural identification of $S_V(r)$ and $S_V(1)$.  
Here we note that the term $E$ in the right-hand side of $(3.1)$ (of [AF]) vainishes in our setting beacuse $E$ 
is defined by $E=\langle w,\nu-e\rangle$ in Page 283 of [AF] and, in our setting, $w$ is equal to ${\bf 0}$ 
by the $W$-invariantness of $\phi_t$.  
According to Lemma 3.6 of [AF], there exists a positive constant $R_0$ such that, if $r<R_0$, 
then the solution $\widehat r_t$ of the evolution equation $(4.1)$ satisfying the initial condition 
$\widehat r_0=r$ uniquely exists in infinite time.  According to the discussion in 
Page 299(Line 3 from the bottom)-300(Line 8) of [AF], $\widehat r_{t_i}$ 
%($t\in[0,\infty)$) remain to be strictly convex and 
converges to a $W$-equivariant $C^{\infty}$-function $\widehat r_{\infty}$ over $S_V(1)$ 
(in the $C^{\infty}$-topology) as $t\to\infty$ for some sequence $\{t_i\}_{i=1}^{\infty}$ in $[0,\infty)$ 
with $\lim_{i\to\infty}t_i=\infty$.  
This fact implies that the solution $\phi_t$ of $(E_{\cal S})$ satisfying $\phi_0=\iota_r$ exists uniquely 
in infinite time and that $\phi_{t_i}$ converges to a $W$-equivariant 
$C^{\infty}$-embedding $\phi_{\infty}$ of $S_V(r)$ into $V$ (in the $C^{\infty}$-topology) as $i\to\infty$.  
The positive constant $\delta_3$ (which corresponds to the above $R_0$) in Lemma 3.6 of [AF] is smaller than 
the positive constant $\delta_2$ in Lemma 3.2 of [AF]
%(see P291 of [AF])
, $\delta_2$ is smaller than the half of the positive constant $\delta_1(=\delta_M)$ in Lemma 1.1 of [AF] 
and $\delta_1$ is smaller than the half of the positive constant $\delta_0$ defined in Page 257 of [AF].  
Also, according to the definition of $\delta_0$ (see P257 of [AF]), we see that the positive constant 
corresponding to $\delta_0$ is equal to $\frac{r_{\cal S}}{2}$ in our setting.  
Thus $R_0$ is smaller than $\frac{r_{\cal S}}{8}$.  
Denote by $P(\widehat h_t)$ the right-hand side of the evolution equation in Lemma 4.2.  
Assume that $v\in{\rm Ker}\,(\widehat h_t)_Z$.  We may assume that $Z\in S_V(r)\cap\overline C$ without loss of 
generality.  Then we have 
$$P(\widehat h_t)_Z(v,v)=(\widehat{\nabla}^td(\rho_{{\cal S},\phi_t}))_Z(v,v).\leqno{(4.3)}$$
Let $\gamma$ be the $\nabla^t$-geodesic in $S_V(r)$ with $\gamma'(0)=v$.  
Then we have 
$$\begin{array}{l}
\displaystyle{(\widehat{\nabla}^td(\rho_{{\cal S},\phi_t}))_Z(v,v)
=\left.\frac{d^2}{ds^2}\right\vert_{s=0}\rho_{{\cal S},\phi_t}(\gamma(s))}\\
\displaystyle{=\sum_{\alpha\in\triangle_+}m_{\alpha}
\left(\frac{\sqrt{\varepsilon}
\alpha(\widetilde{\nabla}^{\nu_t\circ\gamma}_{\frac{\partial}{\partial s}\vert_{s=0}}\nu_{t\ast}(\gamma'(s)))}
{\tan(\sqrt{\varepsilon}\alpha(\phi_t(Z)))}
-\frac{2\varepsilon\alpha(\phi_t(Z))\alpha(\nu_{t\ast}(v))}
{\sin^2(\sqrt{\varepsilon}\alpha(\phi_t(Z)))}\right.}\\
\hspace{0.5truecm}\displaystyle{
-\frac{\sqrt{\varepsilon}^2\alpha(\nu_t(Z))
\alpha(\widetilde{\nabla}^{\phi_t\circ\gamma}_{\frac{\partial}{\partial s}\vert_{s=0}}
\phi_{t\ast}(\gamma'(s)))}{\sin^2(\sqrt{\varepsilon}\alpha(\phi_t(Z)))}}\\
\hspace{0.5truecm}\displaystyle{\left.
+\frac{2\sqrt{\varepsilon}^3\alpha(\phi_{t\ast}(v))^2\alpha(\nu_t(Z))}
{\sin^2(\sqrt{\varepsilon}\alpha(\phi_t(Z)))\tan(\sqrt{\varepsilon}\alpha(\phi_t(Z)))}\right)}\\
%\end{array}$$
%\newpage
%$$\begin{array}{l}
\displaystyle{=\sum_{\alpha\in\triangle_+}m_{\alpha}
\left(\frac{\sqrt{\varepsilon}\alpha(\phi_{t\ast}((\widehat{\nabla}^t_v\widehat A_t)(v)))
-\widehat h_t(\widehat A_t(v),v)\alpha(\nu_t(Z))}
{\tan(\sqrt{\varepsilon}\alpha(\phi_t(Z)))}\right.}\\
\hspace{0.5truecm}\displaystyle{
-\frac{2\varepsilon\alpha(\phi_t(Z))\alpha(\phi_{t\ast}(\widehat A_t(v)))}
{\sin^2(\sqrt{\varepsilon}\alpha(\phi_t(Z)))}
-\frac{\varepsilon\alpha(\nu_t(Z))^2\widehat h_t(v,v)}
{\sin^2(\sqrt{\varepsilon}\alpha(\phi_t(Z)))}}\\
\hspace{0.5truecm}\displaystyle{\left.
+\frac{2\sqrt{\varepsilon}^3\alpha(\phi_{t\ast}(v))^2\alpha(\nu_t(Z))}
{\sin^2(\sqrt{\varepsilon}\alpha(\phi_t(Z)))\tan(\sqrt{\varepsilon}\alpha(\phi_t(Z)))}\right),}
\end{array}$$
where $\widetilde{\nabla}^{\phi_t\circ\gamma}$ (resp. $\widetilde{\nabla}^{\nu_t\circ\gamma}$) 
is the pullback connection of the connection $\widetilde{\nabla}$ of $V$ 
by $\phi_t\circ\gamma$ (resp. $\nu_t\circ\gamma$).  
Hence, since $v\in{\rm Ker}\,\widehat h_t$, we obtain 
$$\begin{array}{l}
\displaystyle{(\nabla^td(\rho_{{\cal S},\phi_t}))_Z(v,v)
=\sum_{\alpha\in\triangle_+}\frac{m_{\alpha}\sqrt{\varepsilon}}
{\tan(\sqrt{\varepsilon}\alpha(\phi_t(Z)))}}\\
\hspace{3.5truecm}\displaystyle{\times\left(
\alpha(\phi_{t\ast}((\widehat{\nabla}^t_v\widehat A_t)(v)))
+\frac{2\varepsilon\alpha(\phi_{t\ast}(v))^2\alpha(\nu_t(Z))}
{\sin^2(\sqrt{\varepsilon}\alpha(\phi_t(Z)))}\right).}
\end{array}\leqno{(4.4)}$$
According to the proof of Lemma 3.6 of [AF], we may assume that 
$$\mathop{\sup}_{t\in[0,\infty)}\vert\vert\vert\phi_t(Z)\vert\vert<R_0,\qquad
\mathop{\sup}_{t\in[0,\infty)}\vert\vert\phi_{t\ast}(\widehat A_t(v))\vert\vert<\varepsilon_1$$
and 
$$\mathop{\sup}_{t\in[0,\infty)}\vert\vert\phi_{t\ast}(v)\vert\vert<\vert\vert v\vert\vert+\varepsilon_2$$
for sufficiently small positive constants $\varepsilon_1$ and $\varepsilon_2$, 
where we note that the statement of Lemma 3.6 is done in the space 
$W^T_{\bullet\bullet}$.  See Page 251-252 of [AF] about the definition of $W^T_{\bullet\bullet}$, 
where we note that $T=\infty$ according to the statement of this lemma.  
By imitating the discussion in the proof of Lemma 8.3 of [Hu1], it follows from the above second inequality that 
$$\mathop{\sup}_{t\in[0,\infty)}\vert\vert\phi_{t\ast}((\widehat{\nabla}^t_v\widehat A_t)(v))\vert\vert
<\varepsilon'_1,$$
where $\varepsilon'_1$ also is a sufficiently small positive constant because $\varepsilon_1$ is sufficiently 
small.  Also, we see that $\alpha(\nu_t(Z))\geq0$ ($t\in[0,\infty)$) because the statement of Lemma 3.6 is done 
in the space $W^{\infty}_{\bullet\bullet}$.  
Hence, by taking $R_0$ as a sufficiently small positive constant, we can show that the right-hand side 
of $(4.4)$ is greater than or equal to zero in all times.  
Hence we have $P(\widehat h_t)_Z(v,v)\geq0\,\,(t\in[0,\infty))$.  
Therefore, since $\phi_0=\iota_r$ is strictly convex (i.e., $\widehat h_0>0$), 
it follows from the maximum principle that 
$\widehat  h_t>0$, that is, $\phi_t$ is strictly convex for all $t\in[0,\infty)$ 
and hence so is also $\phi_{\infty}$.  This completes the proof.  \hspace{11.25truecm}q.e.d.
%\begin{flushright}q.e.d.\end{flushright}

\section{Proof of Theorem B}
In this section, we shall prove Theorem B stated in Introduction.  
We use the notations in Sections 1-4.  
Let ${\cal S}=(V,(\triangle,\{m_{\alpha}\,\vert\,\alpha\in\triangle\},\varepsilon)),\,\phi_t$ and $\iota_r$ 
be in the statement of Theorem B.  

First we shall prove the statement (i) of Theorem B.  

\vspace{0.5truecm}

\noindent
{\it Proof of {\rm(i)} of Theorem B.} 
Let $\phi_t$ ($t\in [0,\infty)$) be the solution of $(E_{\cal S})$ satisfying $\phi_0=\iota_r$.  
The existence of this flow is assured by Theorem A.  
Define a map $f_t$ of the geodesic sphere $S(r):=K\cdot\pi(S_V(r))$ in $G/K$ into $G/K$ by 
$$f_t(k\pi(Z)):=k\pi(\phi_t(Z))\quad(k\in K,\,Z\in S_V(r)).\leqno{(5.1)}$$
%Since $f_t$ is a $K$-equivariant embedding, $\phi_t$ is a $W$-equivariant $C^{\infty}$-embedding.  
Denote by $N_t$ the outward unit normal vector field of $f_t$ and 
$A_t,\,\,H_t,\,\,\overline H_t$ the shape operator, the mean curvature and the average mean curvature of $f_t$ 
for the inward unit normal vector field $-N_t$, respectively.  According to $(3.6)$, we have 
$$(\overline H_t-H_t)N_t\circ\pi\vert_{S_V(r)}=D_{\cal S}(\phi_t),
\leqno{(5.2)}$$
where we use the fact that $\vert\vert\Delta_{\widehat g_t}\phi_t\vert\vert$ is the mean curvature of $\phi_t$.  
Here we note that $(({\overline H}_t-H_t)N_t)_x\in T_{f_t(x)}{\cal T}(=\mathfrak a=V)$ ($x\in\pi(S_V(r))$) 
and hence $(\overline H_t-H_t)N_t\circ\pi\vert_{S_V(r)}$ is regarded as a map from $S_V(r)$ to $V$.  
Since $\phi_t$ ($t\in[0,\infty)$) is the solution of $(E_{\cal S})$ starting from $\iota_r$, it follows from 
$(5.2)$ that $\{f_t\}_{t\in[0,\infty)}$ is the volume-preserving mean curvature flow starting from 
the inclusion map $\iota_{S(r)}:S(r)\hookrightarrow G/K$.  
Hence, since $r<R_0$ by the assumption, 
it follows from Theorem A that $\phi_{t_i}$ converges to a strictly convex embedding $\phi_{\infty}$ 
(in $C^{\infty}$-topology) as $i\to\infty$ for some sequence $\{t_i\}_{i=1}^{\infty}$ in $[0,\infty)$ 
with $\lim_{i\to\infty}t_i=\infty$ and that $\phi_t$ ($0\leq t<\infty$) remain to be strictly convex and hence 
so is also $\phi_{\infty}$.  
Let $f_{\infty}$ be the map the geodesic sphere $S(r)$ into $G/K$ 
defined as in $(5.1)$ for $\phi_{\infty}$ instead of $\phi_t$.  
Since $\phi_{\infty}$ is strictly convex, it follows from $(3.4)$ and $(3.5)$ that so is also $f_{\infty}$.  
Also, since $\{f_t\}_{t\in[0,\infty)}$ is the volume-preserving mean curvature flow, it follows from 
Main theorem of [AF] that $f_{\infty}$ is of constant mean curvature.  
This completes the proof. \hspace{5.75truecm}q.e.d.
%\begin{flushright}q.e.d.\end{flushright}

\vspace{0.5truecm}

%\noindent
%{\it Remark 5.1.} We can show the strictly convexity of $f_t$ (and hence $f_{\infty}$) by using Theorem A 
%%without the use of Main Theorem of [AF] 
%as follows.  It is clear that $f_t(S(r))=K\cdot\pi(\phi_t(S_V(r))$.  
%Since $\phi_t(S_V(r))$ is strictly convex by Theorem A, it follows from $(3.4)$ and $(3.5)$ that 
%so is also $f_t(M)$.  
%\vspace{0.5truecm}

To prove the statements (ii) and (iii) of Theorem B, we prepare the following lemma.  

\vspace{0.5truecm}

\noindent
{\bf Lemma 5.1.} 
{\sl Let $Z^{\max}$ (resp. $Z^{\min}$) be one of maximum (resp. minimum) points of 
$\rho_{{\cal S},\iota_r}$ (hence $H_0\circ\pi\vert_{S_V(r)}$).  
Then the curves 
$t\mapsto\phi_t(Z^{\max})\,\,\,(t\in[0,\infty))$ 
and $t\mapsto\phi_t(Z^{\min})\,\,\,(t\in[0,\infty))$ are described as 
$$\begin{array}{l}
\displaystyle{\phi_t(Z^{\max})=\left(1+\frac{1}{r}\int_0^t(\overline H_t-H_t(\pi(Z^{\max})))dt
\right)Z^{\max}\qquad(t\in[0,\infty)),}\\
\displaystyle{\phi_t(Z^{\min})=\left(1+\frac{1}{r}\int_0^t(\overline H_t-H_t(\pi(Z^{\min})))dt
\right)Z^{\min}\qquad(t\in[0,\infty)).}
\end{array}\leqno{(5.3)}$$
}

\vspace{0.5truecm}

\noindent
{\it Proof.} 
Denote by $\widehat g_t$ and $\nu_t$ the induced metric and the outward unit normal vector field of 
$\phi_t$, respectively.  
First we shall calculate $\displaystyle{\frac{\partial\nu}{\partial t}}$.  
Since $\langle\nu_t,\nu_t\rangle=1$, we have 
$\langle\frac{\partial\nu}{\partial t},\nu_t\rangle=0$.  Hence 
$(\frac{\partial\nu}{\partial t})_Z$ is tangent to 
$\phi_t(S_{V}(r))$ at $\phi_t(Z)$ for each $Z\in S_{V}(r)$.  
Fix $(Z_0,t_0)\in S_{V}(r)\times[0,\infty)$.  
Let $\{e_i\}_{i=1}^{{\it l}-1}$ be an orthonormal base of $T_{Z_0}S_{V}(r)$ 
with respect to $(\widehat g_{t_0})_{Z_0}$ and $\overline e_i$ the tangent vector field of 
$S_{V}(r)\times[0,\infty)$ along $\{Z_0\}\times[0,\infty)$ defined by 
$(\overline e_i)_{(Z_0,t)}:=(e_i)^L_{(Z_0,t)}\,\,\,((Z_0,t)\in\{Z_0\}\times[0,\infty))$, where 
$(e_i)^L_{(Z_0,t)}$ is the horizontal lift of $e_i$ to $(Z_0,t)$ (with respect to the natural 
projection of $S_{V}(r)\times[0,\infty)$ onto $S_{V}(r)$).  
Then we have 
$$\begin{array}{l}
\displaystyle{\left(\frac{\partial \nu}{\partial t}\right)_{(Z_0,t_0)}
=\sum_{i=1}^{{\it l}-1}\left\langle\left(\frac{\partial \nu}{\partial t}\right)_{(Z_0,t_0)},
\phi_{t_0\ast}(e_i)\right\rangle \phi_{t_0\ast}(e_i)}\\
\displaystyle{=-\sum_{i=1}^{{\it l}-1}\left\langle(\nu_{t_0})_{Z_0},
\left(\frac{\partial \phi_{t\ast}(e_i)}{\partial t}\right)_{(Z_0,t_0)}\right\rangle 
\phi_{t_0\ast}(e_i)}\\
\displaystyle{=-\sum_{i=1}^{{\it l}-1}\left\langle(\nu_{t_0})_{Z_0},
\left.\frac{\partial}{\partial t}(\overline e_i \phi)\right\vert_{t=t_0}\right\rangle
\phi_{t_0\ast}(e_i)}\\
\displaystyle{=-\sum_{i=1}^{{\it l}-1}\left\langle(\nu_{t_0})_{Z_0},
e_i\left(\left.\frac{\partial \phi}{\partial t}\right\vert_{t=t_0}\right)\right\rangle 
\phi_{t_0\ast}(e_i)}\\
\displaystyle{=-\sum_{i=1}^{{\it l}-1}\langle(\nu_{t_0})_{Z_0},
e_i(\overline H_{t_0}-(H_{t_0}\circ\pi\vert_{S_{V}(r)}))(\nu_{t_0})_{Z_0}\rangle 
\phi_{t_0\ast}(e_i)}\\
\displaystyle{=\sum_{i=1}^{{\it l}-1}e_i(H_{t_0}\circ\pi\vert_{S_{V}(r)})
\phi_{t_0\ast}(e_i)}\\
\displaystyle{=\sum_{i=1}^{{\it l}-1}\widehat g_{t_0}
(({\rm grad}_{\widehat g_{t_0}}(H_{t_0}\circ\pi\vert_{S_{V}(r)}))_{Z_0},e_i)
\phi_{t_0\ast}(e_i)}\\
\displaystyle{=(\phi_{t_0})_{\ast}
(({\rm grad}_{\widehat g_{t_0}}(H_{t_0}\circ\pi\vert_{S_{V}(r)}))_{Z_0}),}
\end{array}\leqno{(5.4)}$$
where we use $[\frac{\partial}{\partial t},\overline e_i]=0$.  

Now we shall derive $(5.3)$ in terms of $(5.4)$.  
It is clear that there exists a $C^{\infty}$-curve $t\mapsto Z^{\max}_t\,\,\,(t\in[0,\varepsilon))$ 
such that $Z^{\max}_0=Z^{\max}$ and that $Z^{\max}_t$ is a maximum point of 
$H_t\circ\pi\vert_{S_{V}(r)}$ for each $t\in[0,\varepsilon)$, where $\varepsilon$ is a positive 
constant.  Similarly, there exists a $C^{\infty}$-curve 
$t\mapsto Z^{\min}_t\,\,\,(t\in[0,\widehat{\varepsilon}))$ such that $Z^{\min}_0=Z^{\min}$ and that 
$Z^{\min}_t$ is a minimum point of $H_t\circ\pi\vert_{S_{V}(r)}$ for each 
$t\in[0,\widehat{\varepsilon})$, where $\widehat{\varepsilon}$ is a positive constant.  
According to $(5.3)$, for each $t_0\in[0,\varepsilon)$, 
$\frac{\partial\vert\vert \phi_t\vert\vert}{\partial t}\vert_{(Z^{\max}_{t_0},t_0)}<0$ and 
$Z^{\max}_{t_0}$ is a minimum point of $Z\mapsto\frac{d\vert\vert \phi_t\vert\vert}{\partial t}
\vert_{(Z,t_0)}$.  
On the other hand, according to $(5.4)$, we have 
$\frac{d(\nu_t)}{dt}\vert_{Z^{\max}_t}=0$, that is, 
$(\nu_t)_{Z^{\max}_t}=(\nu_0)_{Z^{\max}}$ ($t\in[0,\varepsilon)$).  
From these facts, we can derive that 
$\phi_t(Z^{\max}_t)=\lambda_1(t)Z^{\max}\,\,\,(t\in[0,\varepsilon))$ 
for some positive funcion $\lambda_1$ over $[0,\varepsilon)$ and that 
$Z^{\max}_t=Z^{\max}\,\,\,(t\in[0,\varepsilon))$ holds.  
Similarly, we we can show that 
$\phi_t(Z^{\min}_t)=\lambda_2(t)Z^{\min}\,\,\,(t\in[0,\widehat{\varepsilon}))$ 
for some positive funcion $\lambda_2$ over $[0,\widehat{\varepsilon})$ and that $Z^{\min}_t=Z^{\min}\,\,\,
(t\in[0,\widehat{\varepsilon}))$ holds.   
Hence it follows from $(5.3)$ that 
$$\frac{\partial\lambda_1}{\partial t}
=\frac{1}{r}\left(\overline H_t-H_t(\pi(Z^{\max}))\right).
%\left(\vert\vert(\Delta_{g^{\phi}_t}\phi_t)_{Z^{\max}_0}\vert\vert
%+\rho\left(\lambda_1(t)Z^{\max}_0,\frac{1}{\vert\vert Z^{\max}_0\vert\vert}Z^{\max}_0\right)
%\right)\right).
$$
and 
$$\frac{\partial\lambda_2}{\partial t}
=\frac{1}{r}\left(\overline H_t-H_t(\pi(Z^{\min}))\right).
%\left(\vert\vert(\Delta_{g^{\phi}_t}\phi_t)_{Z^{\min}_0}\vert\vert
%+\rho\left(\lambda_2(t)Z^{\min}_0,\frac{1}{\vert\vert Z^{\min}_0\vert\vert}Z^{\min}_0\right)
%\right)\right).
%\leqno{(4.11)}
$$
Therefore we can derive 
$$\begin{array}{l}
\displaystyle{\phi_t(Z^{\max})=\left(1+\frac{1}{r}\int_0^t(\overline H_t-H_t(\pi(Z^{\max})))dt
\right)Z^{\max}\qquad(t\in[0,\varepsilon))}\\
\displaystyle{\phi_t(Z^{\min})=\left(1+\frac{1}{r}\int_0^t(\overline H_t-H_t(\pi(Z^{\min})))dt
\right)Z^{\min}\qquad(t\in[0,\widehat{\varepsilon})).}
\end{array}
%\leqno{(4.12)}
$$
It is easy to show that these relations hold over $[0,\infty)$.  
This completes the proof.  
%\hspace{10.2truecm}q.e.d.
\begin{flushright}q.e.d.\end{flushright}

\vspace{0.5truecm}

According to this lemma, in the case where $\triangle$ is of type $({\mathfrak a}_2)$, 
the flow $\phi_t(S_{V}(r))$ is as in Figure 3 or 4 for example.  
By using this lemma, we prove the statements (ii) and (iii) of Theorem B.  

\vspace{0.5truecm}

\noindent
{\it Proof of {\rm (ii)} and {\rm (iii)} Theorem B.} 
Let $P$ be a two-plane in $V$ with $P_{\max}\not=\emptyset$ and $P_{\min}\not=\emptyset$.  
Take $Z^{\max}\in P_{\max}$ and $Z^{\min}\in P_{\min,Z^{\max}}$.  
According to Lemma 5.1, $\phi_t(Z^{\max})$ and $\phi_t(Z^{\min})$ are as in $(5.3)$.  
Let $\theta$ be the angle between $\overrightarrow{{\bf 0}Z^{\max}}$ and 
$\overrightarrow{{\bf 0}Z^{\min}}$.  
From the convexity of $\phi_t(S_{V}(r))$, we can derive that 
the part between $\phi_t(Z^{\max})$ and $\phi_t(Z^{\min})$ of the curve 
$\phi_{\infty}(S_{V}(r))\cap P\cap C$ is included by 
$(B_{V}(\frac{r}{\cos\theta})\setminus B_{V}(r\cos\theta))\cap P\cap C$ (see Figure 5).  
From this fact and the definition of $\theta_0$, we can derive 
$M\subset B(\frac{r}{\cos\theta_0})\setminus B(r\cos\theta_0)$.  
Thus the statement (ii) of Theorem B follows.  Since $r$ 
is a sufficiently small positive constant smaller than $R_0$, $\max_{S_V(r)}\widehat r_t$ is sufficiently small 
for all $t\in[0,\infty)$.  Hence we have $3\cos^2(\sqrt{\varepsilon}(\alpha\circ\phi_t))\geq1$ for all 
$t\in[0,\infty)$.  Therefore, according to the maximum principle, it follows from the evolution equation $(4.2)$ 
for $\widehat H^{\cal S}_t$ that 
$\min_{S_V(r)}\,\widehat H^{\cal S}_0\leq \widehat H^{\cal S}_t\leq\max_{S_V(r)}\,\widehat H^{\cal S}_0$, 
which implies that 
$\min_{S(r)}\,H^r\leq H_t\leq\max_{S(r)}\,H^r$ and hence 
$\min_{S(r)}\,H^r\leq H_M\leq\max_{S(r)}\,H^r$.  
On the other hand, according to $(3.8)$, we have $\eta_{\min}(r)\leq H^r\leq\eta_{\max}(r)$.  
Therefore we obtain $\eta_{\min}(r)\leq H_M\leq\eta_{\max}(r)$.  
\hspace{10.4truecm}q.e.d.
%\begin{flushright}q.e.d.\end{flushright}

\vspace{0.4truecm}

\centerline{
%WinTpicVersion3.08
\unitlength 0.1in
%
\hspace{6.5truecm}}

\vspace{1truecm}

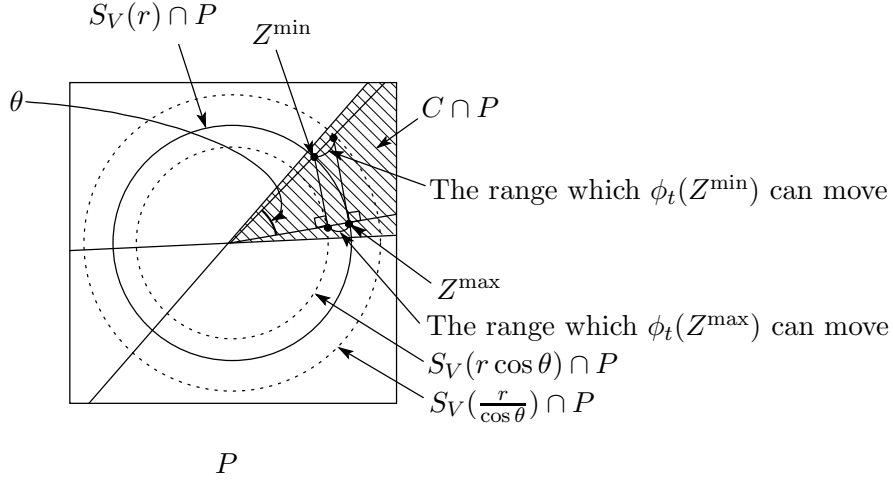
\centerline{{\bf Figure 5.}}

\vspace{0.5truecm}

\noindent
{\it Proof of Corollary C.} 
Let $\theta_0$ be as in the statement (ii) of Theorem B.  
Clearly we have $\theta_0\leq\theta_{G/K}$.  
Hence we obtain the desired inclusion relation from (ii) of Theorem B.  
\hspace{12.5truecm}q.e.d.

\vspace{0.5truecm}

\noindent
{\it Proof of Corollary D.} 
Since $G/K$ is of rank two, we have 
$$\theta_{G/K}=\left\{
\begin{tabular}{ll}
$\displaystyle{\frac{\pi}{3}}$ & (when the root system of $G/K$ is of type $\mathfrak a_2$)\\
$\displaystyle{\frac{\pi}{4}}$ & (when the root system of $G/K$ is of type $\mathfrak b_2$)\\
$\displaystyle{\frac{\pi}{6}}$ & (when the root system of $G/K$ is of type $\mathfrak g_2$).\\
\end{tabular}\right.$$
In the case where $\triangle$ is of type $(\mathfrak a_2)$, it is symmetric for two simple roots 
with considering the multiplicities.  
Hence, in this case, it follows that $\theta_0$ in Theorem B are smaller than or equal to the half 
of $\theta_{G/K}(=$ the half of the length of $S_{V}(1)\cap\overline C$), that is, they are smaller than or equal 
to $\frac{\pi}{6}$ (see Figure 6).  
Therefore, from (ii) of Theorem B, we obtain 
$$M\subset\left\{
%
%
\hspace{8truecm}}

\vspace{1truecm}

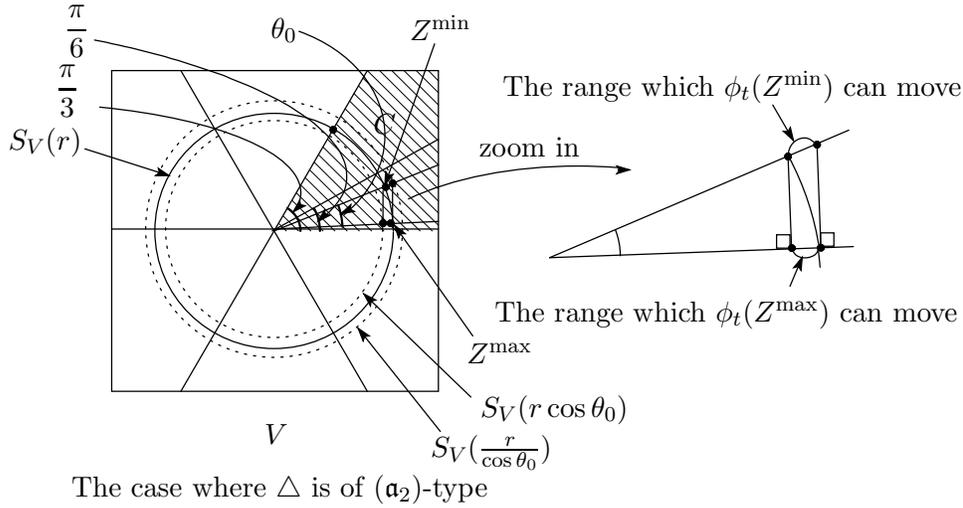
\centerline{{\bf Figure 6.}}

%\newpage

\vspace{0.5truecm}

\centerline{{\bf References}}

\vspace{0.5truecm}

{\small 
\noindent
[AF] N.D. Alikakos and A. Freire, The normalized mean curvature flow for a small bubble 

in a Riemannian manifold, J. Differential Geom. {\bf 64} (2003) 247-303.

\noindent
[CM] E. Cabezas-Rivas and V. Miquel, Volume preserving mean curvature flow in the 

hyperbolic space, Indiana Univ. Math. J. {\bf 56} 
(2007) 2061-2086.

\noindent
[GT] O. Goertsches and G. Thorbergsson, 
On the Geometry of the orbits of Hermann actions, 

Geom. Dedicata {\bf 129} (2007) 101-118.

\noindent
[He] S. Helgason, Differential geometry, Lie groups and symmetric spaces, Academic Press, 

New York, 1978.

\noindent
[Hs] W.Y. Hsiang, On soap bubbles and isoperimetric regions in non-compact symmetric 

spaces, I, Tohoku Math. J. {\bf 44} (1992) 151-175.

\noindent
[HH] W.T. Hsiang and W.Y. Hsiang, On the uniqueness of isoperimetric solutions and 

imbedded soap bubbles in non-compact symmetric spaces, I, Invent. Math. {\bf 98} (1989) 

39-58.

\noindent
[Hu1] G. Huisken, Flow by mean curvature of 
convex surfaces into spheres, J. Differential 

Geom. {\bf 20} (1984) 237-266.

\noindent
[Hu2] G. Huisken, Contracting convex hypersurfaces 
in Riemannian manifolds by their mean 

curvature, Invent. math. {\bf 84} (1986) 463-480.

\noindent
[Hu3] G. Huisken, The volume preserving mean curvature flow, J. reine angew. Math. {\bf 382} 

(1987) 35-48.

\noindent
[K1] N. Koike, Actions of Hermann type and proper complex equifocal submanifolds, 

Osaka J. Math. {\bf 42} (2005) 599-611.

\noindent
[K2] N. Koike, Collapse of the mean curvature flow for equifocal submanifolds, Asian J. 

Math. {\bf 15} (2011) 101-128.

%\noindent
%[M] T. Murphy, Curvature-adapted submanifolds of symmetric spaces, Indiana Univ. Math. 
%
%J. {\bf 61} 831-847 (2012). 

\noindent
[PT] R.S. Palais and C.L. Terng, Critical point theory and submanifold 
geometry, Lecture 

Notes in Math. {\bf 1353}, Springer, Berlin, 1988.

\vspace{0.5truecm}

{\small 
\rightline{Department of Mathematics, Faculty of Science}
\rightline{Tokyo University of Science, 1-3 Kagurazaka}
\rightline{Shinjuku-ku, Tokyo 162-8601 Japan}
\rightline{(koike@ma.kagu.tus.ac.jp)}
}

\end{document}